\documentclass[11pt]{amsart}
\usepackage{amsmath}
\usepackage{amssymb}
\usepackage[usenames]{color}
\usepackage{graphicx}
\usepackage{bm}
\usepackage{enumerate}
\theoremstyle{plain}
\newtheorem{theorem}{Theorem}[section]
\newtheorem{lemma}[theorem]{Lemma}
\newtheorem{prop}[theorem]{Proposition}
\newtheorem{cor}[theorem]{Corollary}

\theoremstyle{definition}
\newtheorem{remark}[theorem]{Remark}
\newtheorem{definition}[theorem]{Definition}

\numberwithin{equation}{section}

\newcommand{\abs}[1]{\left|#1\right|}
\begin{document}




\title[Stability of the Positive Mass Theorem for Axisymmetric Manifolds]{Stability of the Positive Mass Theorem for Axisymmetric Manifolds}

\author{Edward T. Bryden}
\address{Department of Mathematics\\
Stony Brook University\\
Stony Brook, NY 11794, USA}
\email{ebryden@math.sunysb.edu}


\thanks{}

\begin{abstract}
Away from the central axis, we prove the stability of the Positive Mass Theorem in the $W^{1,p}$ sense for asymptotically flat axisymmetric manifolds with nonnegative scalar curvature satisfying some additional technical assumptions. We also derive estimates for the volumes of regions, the areas of axisymmetric surfaces, and the distances between points within the manifolds.
\end{abstract}
\maketitle

\section{Introduction}
\label{sec1} \setcounter{equation}{0}
\setcounter{section}{1}

Based on the formulation of General Relativity, our physical intuition leads us to expect a close relationship between the ADM mass of an asymptotically flat Riemannian manifold and its geometry. Recall that the ADM mass of an asymptotically flat Riemannian manifold is defined to be
\begin{equation}
	m=\lim_{R\rightarrow\infty}\frac{1}{16\pi}\int_{S_R}(g_{ij,j}-g_{jj,i})\nu^{i}.
\end{equation}
In their celebrated Positive Mass Theorem \cite{Scoen-Yau PMT}, Schoen-Yau proved that if an asymptotically flat manifold has nonnegative scalar curvature, then the ADM mass is nonnegative. They
also proved the following rigidity theorem
\begin{equation}
 m=0 \implies M \textrm{ is isometric to Euclidean space.}
\end{equation}
It is natural to ask whether stability also holds; if $M$ has small ADM mass, is $M$ close to
Euclidean space? Lee-Sormani \cite{D Lee-C Sormani ROTSYM} have shown that $M$ need not be smoothly, nor even $C^0$, close to Euclidean space even 
in the spherically symmetric setting; there could be increasingly deep thin gravity wells at the 
center. They conjectured that $M$ is close to Euclidean space in the Sormani-Wenger intrinsic flat (SWIF) sense \cite{Huang Lee Sormani,D Lee-C Sormani ROTSYM}. Proving it will require a method for picking appropriate subregions geometrically and a way to show that these regions converge in the SWIF metric to a subset of Euclidean space.

In \cite{D Lee-C Sormani ROTSYM}, Lee and Sormani study stability in the rotationally symmetric setting. They show that tubular neighborhoods of fixed radius $D$ about coordinate spheres of fixed area $A$ converge to the Euclidean tubular neighborhood of radius $D$ about a sphere of area $A$. Earlier, Lee had proven convergence to Euclidean space outside a compact set in the conformally flat setting \cite{Lee}. Assuming strong conditions on sectional curvature, Corvino has proven that an asymptotically flat manifold with nonnegative scalar curvature and small ADM mass must be diffeomorphic to $\mathbb{R}^3$ \cite{Justin Corvino}. Finster, Bray and Kath have papers bounding the $L^2$ norm of the curvature \cite{Bray Finster,Finster Kath}.  After the Lee-Sormani paper, LeFloch-Sormani \cite{LeFloch Sormani} proved that
metric tensors converge in the $H^1_{loc}$ sense in the rotationally symmetric setting. Huang-Lee-Sormani proved SWIF
convergence in the graph setting and Sormani-Stavrov proved it in the geometrostatic setting. Allen
proved $L^2$ convergence in regions where the Inverse Mean Curvature Flow is smooth \cite{Brian Allen}. 

Here, we will study the question of stability in the presence of axisymmetry. The class of axisymmetric metrics is both flexible enough to model a range of physically interesting phenomena and restricted enough that we have powerful tools at hand that are not available in the most general setting. Recall that the coordinate expression for an axisymmetric metric in cylindrical coordinates is
\begin{equation}\label{cyl-crd-expr-metrics}
g=e^{2\alpha-2u}(d\rho^2+dz^2)+\rho^2e^{-2u}(d\phi+Bd\rho+Adz)^2,
\end{equation}
where all the functions involved depend only on $\rho$ and $z$. The killing field associated with the axisymmetry of $g$ is$\frac{\partial}{\partial\phi}$. Since we will be studying large families of asymptotically flat metrics, it is natural to require that the family satisfy some type of uniform falloff condition.
\begin{definition}\label{defn-unif-asym}
	Let $\mathcal{M}$ be a family of axisymmetric metrics. Suppose we can parameterize $\mathcal{M}$ by the functions $\alpha$, $u$, $A$, and $B$ in cylindrical coordinates (\ref{cyl-crd-expr-metrics}). If there exists constants $C$ and $R_0$ such that if $g$ is a metric in $\mathcal{M}$, then for all $\sqrt{\rho^2+z^2}=r\ge R_0$ we have
	\begin{equation}
	\abs{\partial^{I}u}\le\frac{C}{r^{1+\abs{I}}}
	\end{equation}
	\begin{equation}
	\abs{\partial^{I}\alpha}\le\frac{C}{r^{1+\abs{I}}}
	\end{equation}
	\begin{equation}
	\abs{\partial^{I}A}\le\frac{C}{r^{1+\abs{I}}}
	\end{equation}
	\begin{equation}
	\abs{\partial^{I}B}\le\frac{C}{r^{1+\abs{I}}},
	\end{equation}
	then we shall call $\mathcal{M}$ uniformly asymptotically flat outside of radius $R_0$
\end{definition}
In \cite{Chrusciel}, Chru{\'s}ciel shows that if $(M,g)$ is an axisymmetric manifold which is asymptotically flat, then there are cylindrical coordinates $(\rho,z,\phi)$ in which $g$ takes the form (\ref{cyl-crd-expr-metrics}). 
Suppose that $g$ has the standard asymptotically flat falloff rate:
\begin{equation}
\abs{\partial^{I}(g-\delta_{\mathbb{R}^3})}\le\frac{C}{r^{1+\abs{I}}},
\end{equation}
where $\delta_{\mathbb{R}^3}$ is the Euclidean metric. 
In general the asymptotic falloff of the functions $\alpha$, $u$, $A$, and $B$ will not be as strong as the those given in Definition \ref{defn-unif-asym}. However, we may make an additional assumption on the killing field of $g$ which will imply that the functions $\alpha$, $u$, $A$, and $B$ do have the same falloff as in \ref{defn-unif-asym}. This indicates that there are many families of metrics satisfying the requirements of Definition \ref{defn-unif-asym}.

In Chru{\'s}ciel's construction of cylindrical coordinates, the coordinate functions $\rho$ and $z$ are both solutions to a PDE determined by the metric $g$. Specifically, if we let $\eta$ denote the killing field generating the axisymmetry of $g$ and let $q$ denote the metric on the orbit space induced by $g$, then both $\rho$ and $z$ solve
\begin{equation}\label{rho equation}
	\Delta_q\omega=\Delta_g\omega-\frac{1}{2\abs{\eta}^2_g}<\nabla\omega,\nabla\abs{\eta}^2_g>_g=0.
\end{equation}
In fact, $\rho$ and $z$ are uniquely determined up to conformal maps in the plane. In section two of \cite{Gibbons-Holzegel}, it is noted that if we insist on mapping the axis of symmetry to itself and preserving asymptotic flatness, then $\rho$ is completely fixed. In addition, we can see that $z$ is unique up to translation. This uniqueness justifies our choice to parameterize families of axisymmetric metrics as we did in Definition \ref{defn-unif-asym}

A major obstacle to proving the stability of the Positive Mass Theorem, perhaps the principal one, is that the ADM mass cannot control regions within outermost minimizing surfaces. Classic examples depicting why the Penrose Inequality depends on the area of an outermost minimizing surface demonstrate this phenomenon.
One way to overcome this difficulty, which was 
applied in the work of Bray, Finster, Lee, Kath, Huang-Lee-Sormani, and Allen \cite{Bray Finster,Finster Kath,Huang Lee Sormani,Brian Allen}, is to impose conditions which constrain the location, or prevent the existence, of an outermost minimal surface. We shall follow this approach in making the following definition.
\begin{definition}
	Let $\mathcal{M}$ be a family of axisymmetric metrics and let $\eta$ denote the killing field generating their axisymmetry. Suppose that for each metric $g\in\mathcal{M}$ we have the following inequality
	\begin{equation}\label{defn-area-enl}
	\frac{\abs{\eta}_g}{\abs{\nabla\rho}_g}(\rho_0,z)\ge\rho_0.
	\end{equation}
	Then we shall call $\mathcal{M}$ a family of area enlarging metrics at $\rho_0$. If the inequality holds for each $\rho_0$, then we shall simply call the family area enlarging.
\end{definition}
Uniqueness of solutions to (\ref{rho equation}) implies that the above is a condition imposed on the family $\mathcal{M}$ and has significance beyond a coordinate condition. However, it is useful to express the above in terms of cylindrical coordinates. In coordinates the condition reads
\begin{equation}
	\left(\alpha-2u\right)(\rho_0,z)\ge 0.
\end{equation}

Suppose that $\mathcal{M}$ satisfies condition (\ref{defn-area-enl}) for all $\rho_0$. Let $\delta_{\mathbb{R}^3}$ denote the background Euclidean metric given in the cylindrical coordinates $(\rho,z,\phi)$. Then in Proposition \ref{Area Enlarging} we show that
\begin{equation}
	Area_g(\Sigma)\ge Area_{\delta_{\mathbb{R}^3}}(\Sigma)
\end{equation}
for axisymmetric surfaces $\Sigma$. Together with the Penrose Inequality, the above area inequality works to constrain the location of outermost minimal surfaces. In Corollary \ref{loc outermost min surface} we show that if $\Sigma$ is an axisymmetric outermost minimal surface which is also a sphere, then
\begin{equation}
	\Sigma\subset\rho^{-1}\left([0,2\sqrt{2}m)\right),
\end{equation}
where $m$ is the ADM mass of the metric under consideration.

As in prior work on stability, we must judiciously decide which regions we will study. In view of the above discussion, the regions
\begin{equation}\label{cylindrical region}
	\tilde{\Omega}^{\rho_1}_{\rho_0}(\sigma)=\left\lbrace\rho_0+\sigma\le\rho\le\rho_1,\abs{z}\le\frac{\rho_1}{2}\right\rbrace\times[0,2\pi),
\end{equation}
for some fixed $\rho_0$ and $\sigma\ge0$, are natural choices. If $\sigma$ is identically zero, then we shall write $\tilde{\Omega}^{\rho_1}_{\rho_0}$. Since we mainly work in the orbit space, we shall often only consider the image of $\tilde{\Omega}^{\rho_1}_{\rho_0}(\sigma)$ under the projection map, which is simply the rectangle
\begin{equation}\label{orbit rectangle}
	\Omega^{\rho_1}_{\rho_0}(\sigma)=\left\lbrace\rho_0+\sigma\le\rho\le\rho_1,\abs{z}\le\frac{\rho_1}{2}\right\rbrace.
\end{equation}
Similarly, if $\sigma$ is taken to be zero, then we shall simply write $\Omega^{\rho_1}_{\rho_0}$.

Instead of the area enlarging assumption (\ref{defn-area-enl}), we will at first work with another requirement.
\begin{definition}
	Let $\mathcal{M}$ be a family of axisymmetric metrics. Suppose that for each metric $g\in\mathcal{M}$ we have the following inequality
	\begin{equation}\label{defn-rad-mono}
	\frac{\partial}{\partial\rho}\left(\frac{1}{\rho}\frac{\abs{\eta}_g}{\abs{\nabla\rho}_g}\right)\le 0
	\end{equation}
	on the set $\lbrace\rho=\rho_0\rbrace$.
	Then we shall call the family radially monotone at $\rho_0$. If $\mathcal{M}$ is radially monotone at each $\rho_0$, then we will simply call $\mathcal{M}$ radially monotone.
\end{definition}
This too is a geometric condition on a family of axisymmetric metrics. In  Proposition \ref{Sub-IMCF Radial Monotonicity} we show that if $g$ is an axisymmetric metric,  $\rho$ is the solution to the equation (\ref{rho equation}), then $g$ is radially monotone if and only if the level sets of the function $\rho$ form a sub-inverse-mean-curvature flow.  

The radial monotonicity condition has a useful expression in cylindrical coordinates:
\begin{equation}
	\frac{\partial(\alpha-2u)}{\partial\rho}\le 0.
\end{equation}
In this form, a similar inequality to the above can be found in section 3.2 of \cite{Chrusciel Nguyen}.

One could wonder if there is any relationship between the area enlarging condition and the radial monotonicity condition. The author does not know of any such relationship holding pointwise. However, if radial monotonicity holds everywhere, then the area enlarging condition must also hold everywhere, see Proposition \ref{Radially Monotone implies Area Enlarging}.

In the appendix, we will show that the Kerr-Newman and axisymmetric geometrostatic metrics satisfy radial monotonicity and the area enlarging condition, respectively. In fact, the Kerr-Newman metrics satisfy radially monotonicity strictly, so that small perturbations of the Kerr-Newman metrics are also radially monotone. The same is true for  small perturbations of axisymmetric geometrostatic metrics with regards to the area enlarging condition.  

We now state the stability of the Positive Mass theorem in the $W^{1,p}$ sense.
\begin{theorem}\label{W1p estimates}
	Let $\mathcal{M}$ be a family of axisymmetric metrics with nonnegative scalar curvature which is uniformly asymptotically flat outside of radius $R_0$. Suppose that $\mathcal{M}$ is radially monotone at $\rho_0$ and that for each metric in $\mathcal{M}$, we have
	\begin{equation}
	A=B=0.
	\end{equation}
	For every $\rho_1>\max\lbrace\rho_0,R_0\rbrace$, $\epsilon>0$, $\sigma>0$, and $1\le p<2$ there exists a $\delta>0$ such that if the ADM mass of $g\in\mathcal{M}$ is less than $\delta$, then
	\begin{equation}
	\abs{\abs{g-\delta_{\mathbb{R}^3}}}_{W^{1,p}\left(\tilde{\Omega}^{\rho_1}_{\rho_0}(\sigma)\right)}<\epsilon,
	\end{equation}
	and
	\begin{equation}
		\abs{\abs{q-\delta_{\mathbb{R}^2}}}_{W^{1,p}\left({\Omega}^{\rho_1}_{\rho_0}(\sigma)\right)}<\epsilon,
	\end{equation}
	where $\delta_{\mathbb{R}^3}$ denotes the Euclidean metric in cylindrical coordinates, $\delta_{\mathbb{R}^2}$ denotes the Euclidean metric in the $(\rho,z)$ plane, and $q$ denotes the orbit metric of $g$ in the $(\rho,z)$ plane. $\tilde{\Omega}^{\rho_1}_{\rho_0}(\sigma)$ denotes the cylinder given in (\ref{cylindrical region}) and $\Omega^{\rho_1}_{\rho_0}(\sigma)$ denotes its orbit space.
\end{theorem}
The assumption that the functions $A$ and $B$ vanish is very likely unnecessary, however it does simplify the analysis considerably. That the exponent $p$ is required to be less than two is natural to the problem at hand. Suppose we were able to prove an analogous result for $p>2$. Then, we would be able to apply the Sobolev Embedding Theorem to conclude that the convergence was actually $C_0$ convergence. However, as mentioned before, see \cite{D Lee-C Sormani ROTSYM}, there are counter examples to $C_0$ stability.

It is not yet known if $W^{1,p}$ convergence implies SWIF convergence. However, in the course of proving $W^{1,p}$ stability, we obtain similar estimates to those Huang-Lee-Sormani use to prove the stability of the Positive Mass Theorem in the SWIF metric for graphical manifolds \cite{Huang Lee Sormani}. Let $\mathcal{M}$ be a family of three dimensional asymptotically flat graphical manifolds in $\mathbb{R}^4$ and let $C_{r_0}$ denote the infinite cylinder with base a ball of radius $r_0$ about the origin in $\mathbb{R}^3\subset\mathbb{R}^4$. Huang-Lee-Sormani studied the regions $\Omega_{r_0}\subset M\in\mathcal{M}$ defined by
\begin{equation}
	\Omega_{r_0}:=M\cap C_{r_0},
\end{equation}
for some appropriately large $r_0$. Additionally, they assume a uniform diameter bound on the $\Omega_{r_0}$. They then show that as the ADM mass approaches zero, the regions $\Omega_{r_0}$ converge in the SWIF metric to a three dimensional Euclidean ball in $\mathbb{R}^4$,
\begin{equation}
	B(0,r_0)\times\lbrace 0\rbrace.
\end{equation} 
Their proof follows from three assertions. First, they showed that the volumes of the $\Omega_{r_0}$ converge to the volume of $B(0,r_0)$. Second, they showed that the area of $\partial\Omega_{r_0}$ approaches the area of $\partial B(0,r_0)$. Finally, they showed that $\partial\Omega_{r_0}\cap\partial C_{r_0}$ Lipschitz converges to $\partial B(0,r_0)\times\lbrace 0\rbrace$.  

We are able to establish volume convergence for the cylinders $\tilde{\Omega}^{\rho_1}_{\rho_0}(\sigma)$ defined as in (\ref{cylindrical region}).
\begin{theorem}\label{Volume Estimates}
	Let $\mathcal{M}$ be a family of axisymmetric metrics with nonnegative scalar curvature which is uniformly asymptotically flat outside of radius $R_0$. Suppose also that $\mathcal{M}$ is radially monotone at $\rho_0$.
	For any constants $\epsilon>0$, $\sigma>0$, and $\rho_1>\max\lbrace\rho_0,R_0\rbrace$, 
	there exists a $\delta>0$ such that if $g\in\mathcal{M}$ and
	\begin{equation}
	m(g) < \delta,
	\end{equation} 
	then
	\begin{equation}
	\abs{\Omega}+\epsilon\ge vol_g(\Omega)\ge\abs{\Omega}-\epsilon
	\end{equation}
	for any region $\Omega$ such that
	\begin{equation}
	\Omega \subset \tilde{\Omega}^{\rho_1}_{\rho_0}(\sigma).
	\end{equation}
\end{theorem}

We are also able to establish control over areas inside our designated regions.
\begin{theorem}\label{Area Estimates}
	Let $\mathcal{M}$ be a family of axisymmetric metrics with nonnegative scalar curvature which is  uniformly asymptotically flat outside of radius $R_0$. Suppose also that $\mathcal{M}$ is radially monotone at $\rho_0$.
	For any fixed axisymmetric surface $\Sigma$, constant $\epsilon>0$, and constant $\rho_1>\max\lbrace\rho_0,R_0\rbrace$, there exists a $\delta>0$ such that if $m(g)<\delta$, then
	\begin{equation}
	\abs{\Sigma\cap\tilde{\Omega}^{\rho_1}_{\rho_0}(\sigma)}+\epsilon
	\ge Area_g\left(\Sigma\cap\tilde{\Omega}^{\rho_1}_{\rho_0}(\sigma)\right)
	\ge\abs{\Sigma\cap\tilde{\Omega}^{\rho_1}_{\rho_0}(\sigma)}-\epsilon.
	\end{equation}
\end{theorem}

We obtain an estimate on distances between certain points in $\tilde{\Omega}^{\rho_1}_{\rho_0}(\sigma)$ which can be used to give an upper bound on the diameter of $\tilde{\Omega}^{\rho_1}_{\rho_0}(\sigma)$.
\begin{theorem}\label{d estimates}
Let $\mathcal{M}$ be a family of axisymmetric metrics with nonnegative scalar curvature which is uniformly asymptotically flat outside of radius $R_0$. Suppose $\mathcal{M}$ is also radially monotone at $\rho_0$. Additionally, assume that $A=B=0$ in the coordinate representations of the metrics under consideration. Suppose we are given $\epsilon>0$, $\sigma>0$, and $\rho_1>\max\lbrace\rho_0,R_0\rbrace$. There exists a constant $\delta>0$ such that if $m(g)\le\delta$ and $x$ and $y$ are any points such that the Euclidean line segment connecting them lies in $\Omega^{\rho_1}_{\rho_0}(\sigma)\times\lbrace\phi_0\rbrace$ for any $\phi_0$, then 
\begin{equation}
d_g(x,y)\le d(x,y)+\epsilon
\end{equation}
\end{theorem}
For more general pairs of points $x$ and $y$ in $\tilde{\Omega}^{\rho_1}_{\rho_0}$ we have a pointwise estimate on their distance to each other.
\begin{theorem}\label{d estimates pointwise}
	Let $\mathcal{M}$ be a family of axisymmetric metrics with nonnegative scalar curvature which is uniformly asymptotically flat outside of radius $R_0$. Suppose also that $\mathcal{M}$ is radially monotone at $\rho_0$. Additionally, assume that $A=B=0$ in the coordinate representations of the metrics under consideration. Suppose we are given $\epsilon>0$ and $\sigma>0$ and points $x$ and $y$ such that the Euclidean line segment connecting them lies in $\tilde{\Omega}^{\rho_1}_{\rho_0}(\sigma)$. There exists a constant $\delta>0$ such that if $m(g)\le\delta$, then 
	\begin{equation}
	d_g(x,y)\le d(x,y)+\epsilon
	\end{equation}
\end{theorem}
Finally, we are able to establish uniform convergence at large distances from the origin.
\begin{theorem}\label{ufm con. in the ext.}
	Let $\mathcal{M}$ be a family of axisymmetric metrics with nonnegative scalar curvature which is uniformly asymptotically flat outside of radius $R_0$. Suppose that $\mathcal{M}$ is radially monotone and that for all $g\in\mathcal{M}$ we have 
	\begin{equation}
		A=B=0.
	\end{equation}
	Let $R_1>R_0$ and let $A(R_0,R_1)$ denote the coordinate spherical annulus centered at the origin. For any given $0<\beta<1$ and $\epsilon>0$ there exists a $\delta>0$ such that if $g\in \mathcal{M}$ and
	\begin{equation}
		m(g)<\delta,
	\end{equation}
	then
	\begin{equation}
		\abs{\abs{g-\delta_{\mathbb{R}^3}}}_{C^{0,\beta}\left(A(R_0,R_1)\right)}<\epsilon.
	\end{equation}
\end{theorem}
These theorems are proven in the final section of this paper after we prove a series of Lemmas estimating various terms in the coordinate system. All of the above theorems hold if we assume the area enlarging condition (\ref{defn-area-enl}) instead of radial monotonicity (\ref{defn-rad-mono}). The only change is that in addition to assuming (\ref{defn-area-enl}), we must assume that our family of manifolds satisfies a stronger uniform asymptotic falloff than the one given in Definition \ref{defn-unif-asym}.
\begin{definition}
	Let $\mathcal{M}$ be an uniformly asymptotically flat family of metrics. Suppose that in addition to the uniform asymptotic falloff \ref{defn-unif-asym}, we have some uniform $\tau>0$ such that
	\begin{equation}\label{Str-Unf-Asym-Flat}
	\abs{\alpha}\le\frac{C}{r^{1+\tau}}.
	\end{equation}
	Then we shall call $\mathcal{M}$ strongly uniformly asymptotically flat.
\end{definition}

In the future we would like to prove the Lee-Sormani stability Conjecture that regions outside outermost minimizing surfaces converge in the SWIF sense to regions in Euclidean space.  Our volume, area, and distance controls should be useful towards such a proof. Here we used an extra condition (\ref{defn-area-enl}) to constrain, a priori, the location of outer most minimal surfaces. Another approach would be to actually locate outermost minimal surfaces without any assumption. This was done easily in Lee-Sormani thanks to spherical symmetry and was a huge challenge in the work of Sormani-Stavrov \cite{Sormani Stavrov}. Locating the outermost minimal surfaces in an axisymmetric manifold is of independent interest and would be worthy of a paper on its own.

\medskip
\textbf{Acknowledgments.}
I would like to thank my advisor, Marcus Khuri, for his guidance and for proposing this problem. I'm grateful to Ye Sle Cha, Piotr Chru{\'s}ciel, and Luc Nguyen for their interest in the final result. I'm grateful to Christina Sormani for the support she gives young mathematicians and for the workshops she organizes (funded by DMS 1309360). I'm grateful to Brian Allen, Lisa Hernandez, Sajjad Lakzian, and Dan Lee for the discussions we had at those workshops.


\section[short title]{Background Information}
\label{sec2} \setcounter{equation}{0}
\setcounter{section}{2}
The ADM mass is calculated by taking a limit of integrals over the boundaries of increasingly large coordinate balls. It is thus unclear how the ADM mass should control the geometry inside of these balls. In fact, arbitrary local perturbations of a metric would not change its ADM mass. However, if we restrict our attention to metrics with nonnegative scalar curvature, then we are no longer entirely free in our choice of local perturbation. This restores our hope that the ADM mass can control geometry.

In an attempt to relate ADM mass and the interior geometry, it is natural to make use of the divergence theorem,
\begin{equation}
	m(g)=\lim_{R\rightarrow\infty}\frac{1}{16\pi}\int_{\partial B_R}(g_{ij,j}-g_{jj,i})\nu^{i}=\lim_{R\rightarrow\infty}\frac{1}{16\pi}\int_{B_R}div(g_{ij,j}-g_{jj,i}),	
\end{equation} 
to get an integral over the interior. For now, we are ignoring the question of which metric we should use to take the divergence. Intuitively, we think of scalar curvature as a local energy density. As such, we would like to relate the divergence term to the scalar curvature. Ideally, the nonnegativity of the scalar curvature should give control over the integral of the divergence term. This approach can be successfully carried out in the case of axisymmetric metrics. Furthermore, Witten used a more sophisticated version of this idea to prove the positive mass theorem for manifolds with spinors \cite{Witten Edward}.

In cylindrical coordinates for axisymmetric metrics we have the following formula for the scalar curvature \cite{Brill}:
\begin{equation}\label{Scalar Curvature Formula}
	R_g=4e^{2(u-\alpha)}\left[\Delta_{\mathbb{R}^3}(u-\frac{1}{2}\alpha)-\frac{1}{2}|\nabla u|^2_{\delta}+\frac{1}{2\rho}\frac{\partial\alpha}{\partial\rho}-\frac{\rho^2e^{-2\alpha}}{8}\left(\frac{\partial B}{\partial z}-\frac{\partial A}{\partial\rho}\right)^2\right].
\end{equation}
Here we can see that the scalar curvature is indeed closely related to a divergence, namely $\Delta_{\mathbb{R}^3}(u-\frac{\alpha}{2})$. This observation leads to a very useful formula for the mass \cite{Brill},
\begin{equation}\label{Mass-Formula}
	m(g)=\frac{1}{16\pi}\int_{\mathbb{R}^3}\left[e^{-2(u-\alpha)}\left[R_g+\frac{\rho^2e^{-4\alpha+2u}}{2}\left(\frac{\partial B}{\partial z}-\frac{\partial A}{\partial\rho}\right)^2\right]+2|\nabla u|^2_{\delta}\right]\rho d\rho dzd\phi.
\end{equation}
Since all other terms are explicitly nonnegative, if we assume that $R\ge 0$, then the ADM mass immediately gives control over the gradient of $u$. In an asymptotically flat metric,   $u$ must be arbitrarily small on large coordinate spheres. It is therefore reasonable to suppose that we can use the fundamental theorem of calculus to control $u$ everywhere in the manifold. In order to make this precise, we will use the following representation formula to express $u$ in terms of its gradient and its value on large coordinate spheres.

Suppose $\Omega$ is a compact region on which the divergence theorem holds and let $\Gamma$ be the fundamental solution for the Laplacian. Assume further that $u$ is a function which is differentiable on $CL(\Omega)$. Then we have
\begin{equation}\label{Grdnt-Rep-Frm}
u(x)=-\int_{\partial\Omega}u(y)<\nabla\Gamma(x,y),n>dy+\int_{\Omega}<\nabla u(y),\nabla\Gamma(x,y)>dy.
\end{equation}
In order to see this, we follow the calculations appearing as 2.15 in \cite{Gilbarg-Trudinger}, except we use the divergence theorem on the vector field $Z$ defined by
\begin{equation}
Z=u(y)\nabla\Gamma(x,y).
\end{equation}

The ease with which we can obtain estimates for $u$ is encouraging, however there is one more hurdle. If we want to use mass to control the metric, then we must be able to turn our estimates for $u$ into estimates for $e^u$. Luckily, we may use the well known Moser-Trudinger inequality \cite{Gilbarg-Trudinger} to accomplish this. 

In view of the coordinate expression for an axisymmetric metric (\ref{cyl-crd-expr-metrics}), we know that if we can control $e^{\alpha-2u}$ as well as $e^u$, then we have achieved good control over the metric. Although it is less clear, it is possible to use the mass formula (\ref{Mass-Formula}) and the scalar curvature equation (\ref{Scalar Curvature Formula}) to show that the ADM mass controls the $W^{1,p}$ norm of $\alpha-2u$. The process is similar to what we do to estimate $u$. However, we use Green's representation formula, instead of (\ref{Grdnt-Rep-Frm}), to express $\alpha-2u$ as a boundary term plus an integral of its derivatives. We recall Green's representation formula now. 

Let $\Omega$ be a compact region on which the divergence theorem holds and let $\Gamma$ be the fundamental solution of the Laplacian. Suppose that $\omega$ is a twice differentiable function on $CL(\Omega)$. Then we have the following representation of $\omega$
\begin{equation}\label{Green's Rep-Frm}
	\omega(x)=\int_{\partial\Omega}\left[\omega(y)\frac{\partial\Gamma(x,y)}{\partial\nu}-\Gamma(x,y)\frac{\partial \omega(y)}{\partial \nu}\right]dy+\int_{\Omega}\Gamma(x,y)\Delta \omega(y)dy.
\end{equation}
This result appears in \cite{Gilbarg-Trudinger} as equation 2.16.

With $W^{1,p}$ estimates for $\alpha-2u$ in hand, we might hope to use the Moser-Trudinger inequality to get estimates for $e^{\alpha-2u}$. Unfortunately, the Moser-Trudinger inequality doesn't apply in this case. Luckily, because of axisymmetry, we are essentially working in two dimensions. This gives us extra control that does not exist in higher dimensions. In this setting we are able to prove a result similar to the Moser-Trudinger inequality, which allows us to turn $W^{1,p}$ estimates for $\alpha-2u$ into $W^{1,p}$ estimates for $e^{\alpha-2u}$.

In using (\ref{Grdnt-Rep-Frm}) and (\ref{Green's Rep-Frm}) to control the $W^{1,p}$ norms of $u$ and $\alpha-2u$, we rely on estimates of the Riesz potential. Recall that the Riesz potential of a function $f$ over a region $\Omega$, denoted $(V_\mu f)(x)$, is defined as
\begin{equation}
	(V_\mu f)(x)=\int_{\Omega}\abs{x-y}^{n(\mu-1)}f(y)dy,
\end{equation}
for $\mu\in(0,1]$. Let $0\le\delta=\delta(p,q)=q^{-1}-p^{-1}<\mu$ and let $\omega_n$ denote the volume of the unit $n$ dimensional ball.
The following inequality appears as Lemma 7.12 in \cite{Gilbarg-Trudinger}.
\begin{equation}\label{Riesz Potential Estimate}
	\abs{\abs{(V_\mu f)}}_{p}\le\left(\frac{1-\delta}{\mu-\delta}\right)^{1-\delta}\omega_n^{1-\mu}\abs{\Omega}^{\mu-\delta}\abs{\abs{f}}_{q}.
\end{equation}


\section[short title]{Sobolev Estimates for $u$ and $e^u$}
In this section we will see in greater detail the steps needed to estimate the $W^{1,p}$ norm of $e^u$ using the mass formula (\ref{Mass-Formula}). Our end goal is to produce estimates over the regions $\Omega^{\rho_1}_{\rho_0}(\sigma)$, see (\ref{orbit rectangle}). In fact, we are always able to take $\sigma$ to be zero. To simplify notation, such rectangles will be denoted by $\Omega^{\rho_1}_{\rho_0}$. 

To start, the ADM mass only explicitly bounds the $L^2(\mathbb{R}^3)$ norm of $\nabla u$. The following Lemma demonstrates that this is enough to get $W^{1,2}(B_{r_0})$ control over $u$ for a ball of fixed radius $r_0$ about the origin in $\mathbb{R}^3$.
\begin{lemma}\label{estimate for the Sobolev Norm of u in R3}
	Let $\mathcal{M}$ be a family of axisymmetric metrics with nonnegative scalar curvature which is  uniformly asymptotically flat outside of radius $R_0$, and let $B_{r_0}$ be the ball of radius $r_0$ about the origin. For any $\epsilon>0$ there exists a $\delta>0$ such that if $g\in\mathcal{M}$ and
	\begin{equation}
		m(g) <\delta,
	\end{equation}
	then
	\begin{equation}
	||u||_{W^{1,2}(B_{r_0})}<\epsilon.
	\end{equation}
\end{lemma}

\textbf{Proof}: We note once again that control over $\abs{\abs{\nabla u}}_{L^2(B_{r_0})}$ is an immediate consequence of the mass formula and the nonnegative scalar curvature assumption. In the calculations that follow we will denote the volume of a three dimensional unit ball by $\omega_3$. First, we look at some very large coordinate ball $B(0,r_1)$ with $r_1>\max\lbrace r_0,R_0\rbrace$. If we let $\Gamma$ be the fundamental solution for the Laplacian, then using (\ref{Grdnt-Rep-Frm}) we may express $u$ as
\begin{equation}
u(x)=-\int_{\partial B(0,r_1)}u(y)<\nabla\Gamma(x,y),n>dy+\int_{B(0,r_1)}<\nabla u(y),\nabla\Gamma(x,y)>dy
\end{equation}

Taking the absolute value of both sides and using the triangle inequality on the right hand side shows us that
\begin{equation}\label{key}
|u(x)|\le\int_{\partial B(0,r_1}\frac{|u(y)|}{3\omega_3|x-y|^2}dy +\int_{B(0,r_1)}\frac{|\nabla u(y)|}{3\omega_3|x-y|^2}dy.
\end{equation}
We now integrate $|u|^2$ over $B(0,r_0)$ and use the well known inequality
\begin{equation}
	(a+b)^2\le2(a^2+b^2) \text{ for }a,b\in\mathbb{R}
\end{equation}
to obtain 
\begin{equation}\label{L2 norm of u}
\int_{B(0,r_0)}|u(x)|^2dx\le 2\int_{B(0,r_0)}\left(\int_{\partial B(0,r_1)}\frac{|u(y)|}{3\omega_3|x-y|^2}dy\right)^2+\left(\int_{B(0,r_1)}\frac{|\nabla u(y)|}{3\omega_3|x-y|^2}dy\right)^2dx.
\end{equation}

To bound the second integral on the right hand side we make use of the mass formula (\ref{Mass-Formula}) and the Riesz potential estimate (\ref{Riesz Potential Estimate}) with $\mu=\frac{1}{3}$ and $q=p=2$ to get
\begin{equation}\label{key}
\int_{B(0,r_1)}\left(\int_{B(0,r_1)}\frac{|\nabla u(y)|}{3\omega_3|x-y|^2}dy\right)^2dx\le 8\pi r_1^2m.
\end{equation}
We estimate the first integral on the right using uniform asymptotic flatness \ref{defn-unif-asym} as follows:
\begin{equation}
\int_{B(0,r_0)}\left(\int_{\partial B(0,r_1)}\frac{|u(y)|}{3\omega_3|x-y|^2}dy\right)^2\le\frac{1}{9\omega_3^2}\int_{B(0,r_0)}\left(\int_{\partial B(0,r_1)}\frac{C}{|x-y|^{2}}\frac{1}{r_1}dy\right)^2\le\frac{\omega_3r_0^3C^2r_1^4}{(r_1-r_0)^4r_1^{2}}.
\end{equation}
Substituting the above two inequalities into (\ref{L2 norm of u}), we obtain
\begin{equation}\label{Estimate of L2 norm of u}
\int_{B(0,r_0)}|u(x)|^2dx\le 2\left[\frac{C^2\omega_3 r_0^3r_1^4}{(r_1-r_0)^4r_1^{2}}+8\pi r_1^2m\right]
\end{equation}
If we let $r_1$ grow arbitrarily large, then the first term on the right will become arbitrarily small. We may counter any growth in the second term on the right by choosing the mass to be small enough. $\qed$

The next step is to estimate $e^u$. In order to do that we will apply the Moser-Trudinger inequality to $u$. Let us now recall the exact statement of the Moser-Trudinger inequality. Let $\Omega\subset\mathbb{R}^n$ and $\omega\in W^{1,n}_0(\Omega)$. Then there exists constants $c_1$ and $c_2$ depending only on $n$, such that
\begin{equation}\label{Msr-Trd Inequality}
\int_{\Omega}\exp\left((\frac{\abs{\omega}}{c_1||\nabla \omega||_n})^{\frac{n}{n-1}}\right)\le c_2|\Omega|.
\end{equation}
This inequality appears as theorem 7.15 in \cite{Gilbarg-Trudinger}. Lemma \ref{estimate for the Sobolev Norm of u in R3} gives $W^{1,2}$ control over $u$, so if we want to apply the Moser-Trudinger inequality, we will have to work over two dimensional domains. Luckily, we have the following almost trivial corollary to Lemma \ref{estimate for the Sobolev Norm of u in R3}.
\begin{cor}\label{Sobolev Convergence for u}
	Let $\mathcal{M}$ be a family of axisymmetric metrics with nonnegative scalar curvature which is  uniformly asymptotically flat outside of radius $R_0$. Let $\Omega^{\rho_1}_{\rho_0}$ denote the region
	\begin{equation}
		 \left\lbrace\rho_0\le\rho\le\rho_1,\abs{z}\le\frac{\rho_1}{2}\right\rbrace.
	\end{equation}
	For every $\epsilon>0$, $\rho_0>0$ and $\rho_1>\rho_0$ there exists a $\delta>0$ such that if the ADM mass of $g\in\mathcal{M}$ is less than $\delta$, then
	\begin{equation}
	\abs{\abs{u}}_{W^{1,2}(\Omega^{\rho_1}_{\rho_0})}<\epsilon.
	\end{equation}
\end{cor}
\textbf{Proof}: Consider the region $\tilde{\Omega}^{\rho_1}_{\rho_0}=\Omega^{\rho_1}_{\rho_0}\times[0,2\pi)$. Choose $r_0$ large enough that
\begin{equation}
	\tilde{\Omega}^{\rho_1}_{\rho_0}\subset B_{r_0}.
\end{equation}
In $\Omega^{\rho_1}_{\rho_0}$ we know that $\rho_0\le\rho$. Thus, we may observe that 
\begin{equation}
	\int_{\Omega^{\rho^1}_{\rho_0}}u^2+\abs{\nabla u}^2d\rho dz\le\frac{1}{2\pi\rho_0}\int_{\tilde{\Omega}^{\rho_1}_{\rho_0}}\left[u^2+\abs{\nabla u}^2\right]\rho d\rho dz d\phi\le\frac{1}{2\pi\rho_0}\abs{\abs{u}}^2_{W^{1,2}(B_{r_0})}.
\end{equation}
Now we may apply Lemma \ref{estimate for the Sobolev Norm of u in R3}.$\qed$

We're now in a position to estimate the $W^{1,p}$ norm of $e^u$. For the $L^p$ norm of $e^u$ the proof is an almost direct application of the Moser Trudinger inequality. To estimate the $L^p$ norm of $\nabla e^u=e^u\nabla u$, we use H{\"o}lder's inequality to analyze each term separately For the $e^u$ term we will once again apply the Moser Trudinger inequality. To estimate $\nabla u$ we will rely on Corollary \ref{Sobolev Convergence for u}.
\begin{lemma}\label{Sobolev Convergence for exp(u)}
	Let $\mathcal{M}$ be a flat family of metrics with nonnegative scalar curvature which is uniformly asymptotically outside of radius $R_0$. Let $\Omega^{\rho_1}_{\rho_0}$ denote the region $\lbrace(\rho,z)|\rho_0\le\rho\le \rho_1,\abs{z}\le\frac{\rho_1}{2}\rbrace$. 
	For every $\rho_1>\rho_0>0$, $\epsilon>0$ and $p<2$ there exists a $\delta>0$ such that if the ADM mass of $g\in\mathcal{M}$ is less than $\delta$, then
	\begin{equation}
	\abs{\abs{e^{\abs{u}}-1}}_{W^{1,p}\left(\Omega^{\rho_1}_{\rho_0}\right)}<\epsilon.
	\end{equation}
\end{lemma}
\textbf{Proof}:
Since $g$ is smooth, $u$ is bounded and has bounded derivatives in $\Omega^{\rho_1}_{\rho_0}$, though we have not made any assumption on what these bounds might be. Thus, $e^{\abs{u}}$ is Lipschitz, and so
\begin{equation}
\int_{\Omega^{\rho_1}_{\rho_0}}|\nabla (e^{|u|}-1)|^p=\int_{\Omega^{\rho_1}_{\rho_0}}|\nabla e^{|u|}|^p=\int_{\Omega^{\rho_1}_{\rho_0}}e^{p\abs{u}}\abs{\nabla u}^p.
\end{equation}
Now, we let $r=\frac{2}{p}$ and apply H{\"o}lder's inequality with $r$ to get
\begin{equation}
\int_{\Omega^{\rho_1}_{\rho_0}}e^{p\abs{u}}\abs{\nabla u}^p\le\left(\int_{\Omega^{\rho_1}_{\rho_0}}e^{r'p\abs{u}}\right)^{\frac{1}{r'}}\left(\int_{\Omega^{\rho_1}_{\rho_0}}\abs{\nabla u}^2\right)^{\frac{p}{2}}.
\end{equation}
Let $D(0,r_0)$ denote the two dimensional disk centered about the origin with radius $r_0$. Choose $r_0$ so that $\Omega^{\rho_1}_{\rho_0}\subset D(0,r_0)$. We may extend $u$ to a function $\bar{u}$ in $W^{1,2}_0(D(0,r_0))$, see theorem 4.7 in \cite{Evans Gariepy}. We may choose the extension $\bar{u}$ such that
\begin{equation}
\abs{\abs{\bar{u}}}_{W^{1,2}_0\left(D(0,r_0)\right)}\le K\abs{\abs{u}}_{W^{1,2}(\Omega^{\rho_1}_{\rho_0})},
\end{equation}
where the constant $K$ is independent of the function $u$.
A quick application of the Cauchy-Schwarz inequality gives us the estimate
\begin{equation}
r'p\abs{\bar{u}}\le\frac{1}{4}\left(r'pc_1\abs{\abs{\nabla \bar{u}}}_2\right)^2+\left(\frac{\abs{\bar{u}}}{c_1\abs{\abs{\nabla \bar{u}}}_2}\right)^2,
\end{equation}
where $c_1$ is the constant appearing in (\ref{Msr-Trd Inequality}).
We may now use the Moser-Trudinger inequality (\ref{Msr-Trd Inequality}) to see that
\begin{equation}\
\left(\int_{\Omega^{\rho_1}_{\rho_0}}e^{r'p\abs{u}}\right)^{\frac{1}{r'}}\le\left(\int_{D(0,r_0)}e^{r'p\abs{\bar{u}}}\right)^\frac{1}{r'}\le\exp\left(\frac{1}{4}r'(pc_1\abs{\abs{\nabla \bar{u}}}_n)^2\right)\left(c_2\abs{D(0,r_0)}\right)^{\frac{1}{r'}}.
\end{equation}
When written entirely in terms of $u$, the above inequality becomes
\begin{equation}\label{exp(u) bound}
	\left(\int_{\Omega^{\rho_1}_{\rho_0}}e^{r'p\abs{u}}\right)^{\frac{1}{r'}}\le\exp\left[\frac{r'}{4}\left(Kpc_1\abs{\abs{ u}}_{W^{1,2}(\Omega^{\rho_1}_{\rho_0})}\right)^2\right]\left(c_2|D(0,r_0)|\right)^{\frac{1}{r'}}.
\end{equation}
Combining this with Corollary (\ref{Sobolev Convergence for u}) gives
\begin{equation}
\int_{\Omega^{\rho_1}_{\rho_0}}|\nabla e^{|u|}|^p\le\exp\left[\frac{r'}{4}\left(Kpc_1\abs{\abs{ u}}_{W^{1,2}(\Omega^{\rho_1}_{\rho_0})}\right)^2\right]\left(c_2\left|D(0,r_0)\right|\right)^{\frac{1}{r'}}\left(\frac{4m}{\rho_0}\right)^{\frac{p}{2}}
\end{equation}

Now that we have successfully estimated $\nabla(e^{\abs{u}}-1)$, we turn to estimating $e^{\abs{u}}-1$. We use the expansion of $e^{\abs{u}}$ to get that
\begin{equation}\label{exp(u) overestimate}
\int_{\Omega^{\rho_1}_{\rho_0}}\abs{e^{\abs{u}}-1}^p=\int_{\Omega^{\rho_1}_{\rho_0}}\left(\sum_{1}^{\infty}\frac{\abs{u}^k}{k!}\right)^p
\end{equation}
Factoring out $\abs{u}$ and over estimating the rest shows that the right hand side is bounded above by
\begin{equation}
\int_{\Omega^{\rho_1}_{\rho_0}}\abs{u}^pe^{p\abs{u}}
\end{equation}
Now, we let $r=\frac{2}{p}$ and apply H{\"o}lder's inequality to get
\begin{equation}
\int_{\Omega^{\rho_1}_{\rho_0}}\abs{u}^pe^{p\abs{u}}\le\left(\int_{\Omega^{\rho_1}_{\rho_0}}\abs{u}^2\right)^{\frac{p}{2}}\left(\int_{\Omega^{\rho_1}_{\rho_0}}e^{r'p\abs{u}}\right)^{\frac{1}{r'}}
\end{equation}
Finally, we may once again apply Lemma \ref{Sobolev Convergence for u} and (\ref{exp(u) bound}) to obtain the result. $\qed$


\section[short title]{Sobolev Estimates for $\alpha-2u$ and $e^{\alpha-2u}$}
We must now concentrate on estimating $\alpha-2u$ and $e^{\alpha-2u}$. We will try to imitate as closely as possible the steps which let us successfully estimate $u$ and $e^{u}$. First, we obtain $W^{1,p}$ estimates for $\alpha-2u$ from the mass formula (\ref{Mass-Formula}). Unfortunately, even at this early stage, the process is notably harder than it was for $u$.

In our attempt to estimate the $W^{1,2}$ norm of $u$ we used a representation formula to express $u$ in terms of its values on a large sphere and its gradient in a large ball. Then we used the asymptotic falloff and the mass formula to control these quantities, respectively. This was a relatively simple process because $\abs{\abs{\nabla u}}$ is a term in the mass formula. However, the gradient of $\alpha-2u$ does not appear directly in the mass formula. Rather, it is the Laplacian of $\alpha-2u$ which appears in the mass formula by way of the scalar curvature equation. We will see the precise nature of this relationship in the following lemmas. For now, the important point is that instead of using (\ref{Grdnt-Rep-Frm}) to express $\alpha-2u$, we should use Green's representation (\ref{Green's Rep-Frm}). It is widely known that one may replace the fundamental solution $\Gamma$ in (\ref{Green's Rep-Frm}) with a function $G(x,y)$, the Green's function of the domain, which vanishes on the boundary of the domain. This choice simplifies Green's representation formula significantly. Unfortunately, the explicit formula for $G(x,y)$ can be complicated depending on the domain. Thus, although our representation formula has been simplified, it is difficult to estimate $G(x,y)$. Luckily, we are working over very simple domains, namely the rectangles $\Omega^{\rho_1}_{\rho_0}$. Therefore, a compromise is possible. We may simplify the representation formula for any one side of the rectangle. Specifically, we may choose a "Green's" function which vanishes, or whose normal derivative vanishes, on one side of the rectangle. Since we have the least amount of a priori knowledge about the metric near the axis of symmetry, we will choose to simplify our representation formula on the side nearest the axis of symmetry.

For the rectangle $\Omega^{\rho_1}_{\rho_0}$, let $\bar{x}$ denote the reflection of the point $x$ about the vertical line $\lbrace\rho=\rho_0\rbrace$. We can define the following two functions
\begin{equation}\label{Neumann rep}
H_{N}(x,y)=\frac{1}{2\pi}\log(\abs{x-y})+\frac{1}{2\pi}\log(\abs{\bar{x}-y})
\end{equation}
and
\begin{equation}\label{Dirichlet rep}
H_D(x,y)=\frac{1}{2\pi}\log(\abs{x-y})-\frac{1}{2\pi}\log(\abs{\bar{x}-y}).
\end{equation}
A quick check shows that we may replace $\Gamma$ by either $H_N$ or $H_D$ in (\ref{Green's Rep-Frm}). Furthermore, a calculation shows that
\begin{equation}
\frac{\partial H_{N}(x,y)}{\partial\nu}|_{\partial\Omega^{\rho_1}_{\rho_0}\cap\lbrace\rho=\rho_0\rbrace}=0
\end{equation}
and
\begin{equation}
H_{D}(x,y)|_{\partial\Omega^{\rho_1}_{\rho_0}\cap\lbrace\rho=\rho_0\rbrace}=0.
\end{equation}
Since we will be integrating against the functions $H_N$ and $H_D$ in what follows, and since $H_N$ and $H_D$ are sums of functions of the form $\log\left(\abs{x-y}\right)$,
it will be useful in what follows to have an $L^p$ estimate for $\log\left(\abs{x-y}\right)$ over bounded regions.
\begin{lemma}
	Let $\Omega$ be a bounded region in $\mathbb{R}^2$ and let
	\begin{equation}
		r_0=\max\lbrace diam(\Omega),1\rbrace.
	\end{equation}
	Then for $y\in cl\left(\Omega\right)$ we have
	\begin{equation}\label{log estimate}
		\int_{\Omega}\abs{\log\left(\abs{x-y}\right)}^{k}dx\le\frac{\pi k!}{2^k}+2\pi(r_0-1)r_0\log(r_0)^k
	\end{equation}
for positive integers $k$.
\end{lemma}
\textbf{Proof:}
We observe that
\begin{equation}
\int_{\Omega}\abs{\log(\abs{x-y})}^kdx\le\int_{B(y,r_0)}\abs{\log(\abs{x-y})}^kdx=\int_{0}^{1}(-1)^k2\pi r\log(r)^kdr+\int_{1}^{r_0}2\pi r\log(r)^kdr
\end{equation}
The second term on the right has the simple estimate
\begin{equation}
	2\pi(r_0-1)r_0\log(r_0)^k.
\end{equation}  
To estimate the first term, one must carry out the integration. By induction, we have the following result.
\begin{equation}
\int_{0}^{1}(-1)^k2\pi r\log(r)^kdr=\frac{\pi k!}{2^k}. \qed
\end{equation}

With all of this in mind, we begin the process of estimating the $W^{1,p}$ norm of $\alpha-2u$.
\begin{prop}\label{Sobolev Convergence of 2u-alpha}
	Let $\mathcal{M}$ be a family of axisymmetric metrics with nonnegative scalar curvature which is  uniformly asymptotically flat outside of radius $R_0$. Suppose in addition that $\mathcal{M}$ is radially monotone at $\rho_0$.
	For every $\rho_1>\rho_0$, $\epsilon>0$ and $p<2$ there exists a $\delta>0$ such that if the ADM mass of $g\in\mathcal{M}$ is less than $\delta$, then
	\begin{equation}
	\abs{\abs{\alpha-2u}}_{W^{1,p}\left(\Omega^{\rho_1}_{\rho_0}\right)}<\epsilon
	\end{equation}
\end{prop}
Applying Green's representation formula to $\alpha-2u$ over the domain $\Omega^{\rho_1}_{\rho_0}$ gives us
\begin{equation}
(\alpha-2u)(x)=\int_{\partial\Omega^{\rho_1}_{\rho_0}}\left[(\alpha-2u)\frac{\partial H_N(x,y)}{\partial\nu}-H_N(x,y)\frac{\partial (\alpha-2u)}{\partial \nu}\right]dy+\int_{\Omega^{\rho_1}_{\rho_0}}H_N(x,y)\Delta(\alpha-2u)dy.	
\end{equation}
The above representation breaks our problem into two pieces. First we must estimate $\Delta(\alpha-2u)$ over $\Omega^{\rho_1}_{\rho_0}$ and then we must estimate $\alpha-2u$ on the boundary of $\Omega^{\rho_1}_{\rho_0}$. The necessary estimates are the content of the following two lemmas.
\begin{lemma}\label{L1 Estimate for the Laplacian of 2u-alpha}
	Let $\mathcal{M}$ be a family of axisymmetric metrics with nonnegative scalar curvature which is uniformly asymptotically flat outside of radius $R_0$. If $g$ is a metric in $\mathcal{M}$ and
	\begin{equation}
	m(g)\le m,
	\end{equation}
	then
	\begin{equation}
	\abs{\abs{\Delta(2u-\alpha)}}_{L^1(\Omega^{\rho_1}_{\rho_0})}\le\frac{4m}{\rho_0}+\frac{4\sqrt{\rho_1 m}}{\rho_0}
	\end{equation}
	for any $\rho_1>\rho_0>0$.
\end{lemma}
\textbf{Proof}: We must relate $\Delta(\alpha-2u)$ to the mass formula. First, we recall that the scalar curvature equation is
\begin{equation}
R_g=4e^{2(u-\alpha)}\left[\Delta_{\mathbb{R}^3}(u-\frac{1}{2}\alpha)-\frac{1}{2}|\nabla u|^2_{\delta}+\frac{1}{2\rho}\frac{\partial\alpha}{\partial\rho}-\frac{\rho^2e^{-2\alpha}}{8}\left(\frac{\partial B}{\partial z}-\frac{\partial A}{\partial\rho}\right)^2\right]
\end{equation}
where we have written $\Delta_{\mathbb{R}^3}$ to emphasize the fact that it is the three dimensional Laplacian which appears, and not the two dimensional Laplacian $\Delta$. However, if we remember that all of the functions involved don't depend on $\phi$, then we can see that
\begin{equation}
	\Delta_{\mathbb{R}^3}(u-\frac{\alpha}{2})=\Delta(u-\frac{\alpha}{2})+\frac{1}{2\rho}\frac{\partial(2u-\alpha)}{\partial\rho}.
\end{equation}
By plugging the above into the scalar curvature equation, we get
\begin{equation}
R_g=4e^{2(u-\alpha)}\left[\Delta(u-\frac{1}{2}\alpha)-\frac{1}{2}|\nabla u|^2_{\delta}+\frac{1}{\rho}\frac{\partial u}{\partial\rho}-\frac{\rho^2e^{-2\alpha}}{8}\left(\frac{\partial B}{\partial z}-\frac{\partial A}{\partial\rho}\right)^2\right].
\end{equation}
We now solve the scalar curvature equation for $\Delta(\alpha-2u)$ and integrate in order to arrive at
\begin{equation}\label{L1 estimate for the laplacian of 2u-alpha}
\int_{\Omega^{\rho_1}_{\rho_0}}\abs{\Delta(\alpha-2u)}d\rho dz\le\int_{\Omega^{\rho_1}_{\rho_0}}\frac{e^{2(\alpha-u)}}{2}R_g+|\nabla u|_{\delta}^2+\frac{2}{\rho}\abs{\frac{\partial u}{\partial\rho}}+\frac{\rho^2e^{-2\alpha}}{4}\left(\frac{\partial B}{\partial z}-\frac{\partial A}{\partial\rho}\right)^2d\rho dz.
\end{equation}
Now, since we are integrating over a region in which $\rho\ge\rho_0$, we have from the mass formula (\ref{Mass-Formula}) that
\begin{equation}
\int_{\Omega^{\rho_1}_{\rho_0}}\frac{e^{2(\alpha-u)}}{2}R_g+|\nabla u|_{\delta}^2+\frac{\rho^2e^{-2\alpha}}{4}\left(\frac{\partial B}{\partial z}-\frac{\partial A}{\partial\rho}\right)^2d\rho dz\le\frac{4m}{\rho_0}.
\end{equation}

To estimate the final term on the right hand side of (\ref{L1 estimate for the laplacian of 2u-alpha}) requires only a little more work. Namely, if we apply H{\"o}lder's inequality to
\begin{equation}
\int_{\Omega^{\rho_1}_{\rho_0}}\frac{2}{\rho}\abs{\frac{\partial u}{\partial\rho}}d\rho dz
\end{equation}
and make the simple estimate $\abs{\frac{\partial u}{\partial\rho}}\le\abs{\nabla u}_{\delta}$, then we obtain
\begin{equation}
\int_{\Omega^{\rho_1}_{\rho_0}}\frac{2}{\rho}\abs{\frac{\partial u}{\partial\rho}}d\rho dz\le\left(\int_{\Omega^{\rho_1}_{\rho_0}}\frac{4}{\rho^2}\right)^{\frac{1}{2}}\left(\int_{\Omega^{\rho_1}_{\rho_0}}\abs{\nabla u}^2_{\delta}d\rho dz\right)^{\frac{1}{2}}.
\end{equation}
Using the mass formula once more, we see that 
\begin{equation}
\left(\int_{\Omega^{\rho_1}_{\rho_0}}\frac{4}{\rho^2}\right)^{\frac{1}{2}}\left(\int_{\Omega^{\rho_1}_{\rho_0}}\abs{\nabla u}^2_{\delta}d\rho dz\right)^{\frac{1}{2}}\le\frac{4\sqrt{\rho_1m}}{\rho_0}.
\end{equation}
Putting each of these estimates together gives the desired result.$\qed$

We now want to estimate boundary terms on $\partial\Omega^{\rho_1}_{\rho_0}$. Due to the asymptotic falloff conditions (\ref{defn-unif-asym}), it is relatively straight forward to estimate terms on $(\partial\Omega^{\rho_1}_{\rho_0})-\lbrace\rho=\rho_0\rbrace$. It is more difficult to estimate terms on $(\partial\Omega^{\rho^1}_{\rho_0})\cap\lbrace\rho=\rho_0\rbrace$.
\begin{lemma}\label{Inner Bdry Est for alpha-2u}
	Let $\mathcal{M}$ be a family of axisymmetric metrics with nonnegative scalar curvature which is uniformly asymptotically flat outside of radius $R_0$. Assume that $\mathcal{M}$ is also radially monotone at $\rho_0$. For $\rho_1>\max\lbrace\rho_0,R_0\rbrace$, if $g\in\mathcal{M}$ and
	\begin{equation}
		m(g)\le m,
	\end{equation}
	then
	\begin{equation}
		\int_{(\partial\Omega^{\rho_1}_{\rho_0})\cap\lbrace\rho=\rho_0\rbrace}\abs{\frac{\partial(\alpha-2u)}{\partial\nu}}\le\frac{4m}{\rho_0}+\frac{4\sqrt{\rho_1m}}{\rho_0}+\frac{6\pi C}{\rho_1},
	\end{equation}
	where the constant $C$ is the one appearing in Definition \ref{defn-unif-asym}.
\end{lemma}
\textbf{Proof}:
It is an easy observation that
\begin{equation}\label{normal derivative to a rectangle}
\frac{\partial}{\partial\nu}|_{\partial\Omega^{\rho_1}_{\rho_0}\cap\lbrace\rho=\rho_0\rbrace}=-\frac{\partial}{\partial\rho}.
\end{equation}
If we write the radial monotonicity condition entirely in terms of coordinate functions, then we may see that for $g\in\mathcal{M}$
\begin{equation}
\frac{\partial(\alpha-2u)}{\partial\rho}(\rho_0,z)\le 0.
\end{equation}
Thus, we observe that
\begin{equation}
\int_{\partial\Omega^{\rho_1}_{\rho_0}\cap\lbrace\rho=\rho_0\rbrace}\abs{\frac{\partial(\alpha-2u)}{\partial\nu}}=-\int_{-\frac{\rho_1}{2}}^{\frac{\rho_1}{2}}\frac{\partial(\alpha-2u)}{\partial\rho}(\rho,z)dz.
\end{equation}
A quick application of Stokes' Theorem over the region 
\begin{equation}
	\lbrace\rho_0\le\rho,\abs{z}\le\frac{\rho_1}{2}\rbrace
\end{equation}
gives
\begin{equation}
\int_{-\frac{\rho_1}{2}}^{\frac{\rho_1}{2}}\frac{\partial(\alpha-2u)}{\partial\rho}=-\int_{\lbrace \rho_0\le\rho,\abs{z}\le\frac{\rho_1}{2}\rbrace}\Delta(\alpha-2u)d\rho dz+\int_{\lbrace\rho\ge\rho_0,\abs{z}=\frac{\rho_1}{2}\rbrace}\frac{\partial(\alpha-2u)}{\partial z}.
\end{equation}
We may estimate the second integral on the right by plugging in the asymptotic estimates (\ref{defn-unif-asym}). The result is the following inequality
\begin{equation}
\abs{\int_{\lbrace\rho\ge\rho_0,\abs{z}=\frac{\rho_1}{2}\rbrace}\frac{\partial(\alpha-2u)}{\partial z}}\le\int_{\lbrace\rho\ge\rho_0,\abs{z}=\frac{\rho_1}{2}\rbrace}\frac{3C}{\abs{(\rho,z)}^2}d\rho.
\end{equation}
We may see by a straightforward integration that
\begin{equation}
\abs{\int_{\lbrace\rho\ge\rho_0,\abs{z}=\frac{\rho_1}{2}\rbrace}\frac{\partial(\alpha-2u)}{\partial z}}\le\frac{6\pi C}{\rho_1}.	
\end{equation}
The last piece of the puzzle is the term
\begin{equation}\label{Laplacian Estimate for the boundary}
\abs{\int_{\lbrace\rho_0\le\rho,\abs{z}\le\frac{\rho_1}{2}\rbrace}\Delta(\alpha-2u)d\rho dz}\le\int_{\lbrace\rho_0\le\rho,\abs{z}\le\frac{\rho_1}{2}\rbrace}|\Delta(\alpha-2u)|d\rho dz.
\end{equation}
We now use the proof of Lemma \ref{L1 Estimate for the Laplacian of 2u-alpha} to bound this term.
Putting everything together, we get
\begin{equation}
\int_{\partial\Omega^{\rho_1}_{\rho_0}\cap\lbrace\rho=\rho_0\rbrace}\abs{\frac{\partial(\alpha-2u)}{\partial\nu}}\le\frac{4m}{\rho_0}+\frac{4\sqrt{\rho_1m}}{\rho_0}+\frac{6\pi C}{\rho_1}.\qed
\end{equation}
We have the necessary estimates to obtain $W^{1,p}$ control over $\alpha-2u$.

\textbf{Proof of Proposition \ref{Sobolev Convergence of 2u-alpha}}: Consider $\Omega^{\tilde{\rho}_1}_{\rho_0}$ for some $\tilde{\rho}_1\ge R_0$. We also choose $\tilde{\rho}_1$ to be much larger than $\rho_1$. As before, we let
\begin{equation}
H_N(x,y)=\frac{1}{2\pi}\log(|x-y|)+\frac{1}{2\pi}\log(|\bar{x}-y|),
\end{equation}
where $\bar{x}$ is the reflection of $x$ about the line $\lbrace\rho=\rho_0\rbrace$. Recall that Green's representation gives us the following formula for $\alpha-2u$:
\begin{equation}
(\alpha-2u)(x)=\int_{\partial\Omega^{\tilde{\rho}_1}_{\rho_0}} (\alpha-2u)(y)\frac{\partial H_N}{\partial\nu}(x,y)-H_N(x,y)\frac{\partial(\alpha-2u)}{\partial\nu}(y)dy+\int_{\Omega^{\tilde{\rho}_1}_{\rho_0}} H_N(x,y)\Delta(\alpha-2u)(y)dy.
\end{equation}
We will imitate the estimates that we made for $u$ in (\ref{Sobolev Convergence for u}). Namely, we see that
\begin{equation}
\int_{\Omega^{\rho_1}_{\rho_0}}\abs{(\alpha-2u)(x)}^{p}dx
\end{equation}
is bounded above by
\begin{equation}
C(p)\int_{\Omega^{\rho_1}_{\rho_0}}\left(\int_{\partial\Omega^{\tilde{\rho}_1}_{\rho_0}}\abs{ H_N\frac{\partial(\alpha-2u)}{\partial\nu}}+\abs{(\alpha-2u)\frac{\partial H_N}{\partial\nu}}dy\right)^p+\left(\int_{\Omega^{\tilde{\rho}_1}_{\rho_0}}\abs{ H_N\Delta(\alpha-2u)}dy\right)^pdx,
\end{equation}
for some constant $C(p)$ depending only on $p$. We estimate each of the three terms above in turn. For the first two terms, we will break $\partial\Omega^{\tilde{\rho}_1}_{\rho_0}$ into
\begin{equation}\label{Asymptotically flat part of the boundary}
\partial\Omega^{\tilde{\rho}_1}_{\rho_0}-\lbrace\rho=\rho_0\rbrace
\end{equation}
and
\begin{equation}\label{Interior part of the boundary}
(\partial\Omega^{\tilde{\rho}_1}_{\rho_0})\cap\lbrace\rho=\rho_0\rbrace.
\end{equation}

Let's start with (\ref{Asymptotically flat part of the boundary}). For this piece of the boundary we can use the uniform asymptotically flat condition to obtain the required estimates. First, notice that for $x$ in $\Omega^{\rho_1}_{\rho_0}$ and $y$ in (\ref{Asymptotically flat part of the boundary}) we have
\begin{equation}
	\abs{H_N(x,y)}\le\frac{\log\left(2diam\left(\Omega^{\tilde{\rho}_1}_{\rho_0}\right)\right)}{\pi}\le\frac{\log(2\sqrt{2}\tilde{\rho_1})}{\pi},
\end{equation}
since $\tilde{\rho}_1$ is much larger than $\rho_0$.
From the asymptotic falloff given in Definition \ref{defn-unif-asym}, we see that for $y$ in (\ref{Asymptotically flat part of the boundary})	
\begin{equation}
	\abs{\frac{\partial(\alpha-2u)}{\partial\nu}(y)}\le\frac{3C}{\tilde{\rho}^2_1}.
\end{equation}
Thus, we may see that
\begin{equation}
\int_{\Omega^{\rho_1}_{\rho_0}}\left(\int_{\partial\Omega^{\tilde{\rho}_1}_{\rho_0}-\lbrace\rho=\rho_0\rbrace}\abs{ H_N(x,y)\frac{\partial(\alpha-2u)}{\partial\nu}(y)}dy\right)^pdx\le\int_{\Omega^{\rho_1}_{\rho_0}}\left(\frac{9\log(2\sqrt{2}\tilde{\rho}_1)C}{\pi\tilde{\rho}_1}\right)^pdx\le \rho_1^2\left(\frac{3\log(2\sqrt{2}\tilde{\rho}_1)C}{\tilde{\rho}_1}\right)^p.
\end{equation}
The other term has a similar estimate:
\begin{equation}
\int_{\Omega^{\rho_1}_{\rho_0}}\left(\int_{\partial\Omega^{\tilde{\rho}_1}_{\rho_0}-\lbrace\rho=\rho_0\rbrace}\abs{(\alpha-2u)\frac{\partial H_N}{\partial\nu}}dy\right)^pdx\le \rho_1^2\left(\frac{6C}{\tilde{\rho}_1-\rho_1}\right)^p.
\end{equation}
Using the two estimates above, we see that
\begin{equation}\label{total asym flat part of boundary estimate}
	\int_{\Omega^{\rho_1}_{\rho_0}}\left(\int_{\partial\Omega^{\tilde{\rho}_1}_{\rho_0}-\lbrace\rho=\rho_0\rbrace}\abs{ H_N\frac{\partial(\alpha-2u)}{\partial\nu}}+\abs{(\alpha-2u)\frac{\partial H_N}{\partial\nu}}dy\right)^p\le C(p)\left(\rho_1^2\left(\frac{3\log(2\sqrt{2}\tilde{\rho}_1)C}{\tilde{\rho}_1}\right)^p+\rho_1^2\left(\frac{6C}{\tilde{\rho}_1-\rho_1}\right)^p\right).
\end{equation}

We can now move to the inner piece of the boundary, (\ref{Interior part of the boundary}). We will further divide $\partial\Omega^{\tilde{\rho}_1}_{\rho_0}\cap\lbrace\rho=\rho_0\rbrace$ into 
\begin{equation}\label{mu bottom}
	\partial\Omega^{\tilde{\rho_1}}_{\rho_0}\cap\lbrace\rho=\rho_0,\abs{z}\le\rho_1\rbrace
\end{equation}
and
\begin{equation}\label{mu top}
	\partial\Omega^{\tilde{\rho}_1}_{\rho_0}\cap\lbrace\rho=\rho_0,\abs{z}\ge\rho_1\rbrace.
\end{equation}
We now estimate
\begin{equation}\label{piece 1 inner boundary}
	\left(\int_{\Omega^{\rho_1}_{\rho_0}}\left(\int_{\partial\Omega^{\tilde{\rho}_1}_{\rho_0}\cap\lbrace\rho=\rho_0,\abs{z}\le\rho_1\rbrace}| H_N(x,y)\frac{\partial(\alpha-2u)}{\partial\nu}|dy\right)^pdx\right)^{\frac{1}{p}}.
\end{equation}
Here we apply Minkowski's inequality for integrals \cite{Folland} to bound the above by
\begin{equation}
\int_{\partial\Omega^{\tilde{\rho}_1}_{\rho_0}\cap\lbrace\rho=\rho_0,\abs{z}\le\rho_1\rbrace}\left(\int_{\Omega^{\rho_1}_{\rho_0}}| H_N(x,y)\frac{\partial(\alpha-2u)}{\partial\nu}|^pdx\right)^{\frac{1}{p}}dy.
\end{equation}
We may rewrite this expression as
\begin{equation}\label{Inner Boundary Rep. Estimate }
	\int_{\partial\Omega^{\tilde{\rho}_1}_{\rho_0}\cap\lbrace\rho=\rho_0,\abs{z}\le\rho_1\rbrace}\abs{\frac{\partial(\alpha-2u)}{\partial\nu}}\left(\int_{\Omega^{\rho_1}_{\rho_0}}|H_N(x,y)|^pdx\right)^{\frac{1}{p}}dy.
\end{equation}
In view of Lemma \ref{Inner Bdry Est for alpha-2u}, we must estimate the $L^P$ norm of $H_N(x,y)$ as a function of $x$ over $\Omega^{\rho_1}_{\rho_0}$ for each $y$ in
\begin{equation}
	\partial\Omega^{\tilde{\rho}_1}_{\rho_0}\cap\lbrace\rho=\rho_0,\abs{z}\le\rho_1\rbrace.
\end{equation}
We see that the points $x$ and $y$ are both contained in $\Omega^{2\rho_1}_{\rho_0}$, Which has diameter $2\sqrt{2}\rho_1$. Let
\begin{equation}
	F(x)=\bar{x}.
\end{equation}
Since $F$ is an isometry, if we apply the change of variable formula to $F$ and note that $y=\bar{y}$ for $y$ in $\lbrace\rho=\rho_0\rbrace$, then we may see that for any $q$, we have
\begin{equation}
	\int_{\Omega^{2\rho_1}_{\rho_0}}\abs{\log\left(\abs{\bar{x}-y}\right)}^qdx=\int_{F\left(\Omega^{2\rho_1}_{\rho_0}\right)}\abs{\log\left(\abs{x-y}\right)}^qdx.
\end{equation} 
Thus, we may use (\ref{log estimate}) to see that
\begin{equation}
\int_{\Omega^{\rho_1}_{\rho_0}}\abs{ H_N(x,y)}dx\le\int_{\Omega^{2\rho_1}_{\rho_0}}\abs{ H_N(x,y)}\le\frac{1}{2}+16\rho_1^2\log(2\sqrt{2}\rho_1),
\end{equation}
and
\begin{equation}
	\left(\int_{\Omega^{\rho_1}_{\rho_0}}\abs{ H_N(x,y)}^2dx\right)^{\frac{1}{2}}\le\left(\int_{\Omega^{2\rho_1}_{\rho_0}}\abs{ H_N(x,y)}^2dx\right)^{\frac{1}{2}}\le\frac{1}{2\pi}\sqrt{2\pi+64\pi\rho_1^2\log(2\sqrt{2}\rho_1)^2}.
\end{equation}
We do a simple interpolation between the above two estimates to get
\begin{equation}
\left(\int_{\Omega^{\rho_1}_{\rho_0}}\abs{ H_N(x,y)}^pdx\right)^{\frac{1}{p}}\le\left(\frac{1}{2}+16\rho_1^2\log(2\sqrt{2}\rho_1)\right)^{\frac{2-p}{p}}\left(\frac{1}{2\pi}\sqrt{2\pi+64\pi\rho_1^2\log(2\sqrt{2}\rho_1)^2}\right)^{\frac{2p-2}{p}},
\end{equation}
We now combine the above with Lemma \ref{Inner Bdry Est for alpha-2u} to bound (\ref{Inner Boundary Rep. Estimate }) by
\begin{equation}\label{Bound on mu bottom}
\left[\frac{1}{2}+16\rho_1^2\log(2\sqrt{2}\rho_1)\right]^{\frac{2-p}{p}}\left[\frac{1}{2\pi}\sqrt{2\pi+64\pi\rho_1^2\log(2\sqrt{2}\rho_1)^2}\right]^{\frac{2p-2}{p}}\left(\frac{4m}{\rho_0}+\frac{4\sqrt{\tilde{\rho}_1m}}{\rho_0}+\frac{6\pi C}{\tilde{\rho}_1}\right)
\end{equation}

The term
\begin{equation}\label{piece 2 inner boundary}
	\left(\int_{\Omega^{\rho_1}_{\rho_0}}\left(\int_{\partial\Omega^{\tilde{\rho}_1}_{\rho_0}\cap\lbrace\rho=\rho_0,\abs{z}\ge\rho_1\rbrace}| H_N(x,y)\frac{\partial(\alpha-2u)}{\partial\nu}|dy\right)^pdx\right)^{\frac{1}{p}}
\end{equation}
is much easier to estimate. In fact, for $x$ in $\Omega^{\rho_1}_{\rho_0}$ and $y$ in $\partial\Omega^{\tilde{\rho}_1}_{\rho_0}\cap\lbrace\rho=\rho_0,\abs{z}\ge\rho_1\rbrace$, we have
\begin{equation}
	\abs{H_N(x,y)}\le \frac{1}{\pi}\max\lbrace\abs{\log(\frac{\rho_1}{2})},\abs{\log(2\sqrt{2}\tilde{\rho}_1)}\rbrace.
\end{equation}
Once again, combining the above with Lemma \ref{Inner Bdry Est for alpha-2u} bounds (\ref{piece 2 inner boundary}) by
\begin{equation}\label{Bound on mu top}
	\frac{(\rho_1)^{\frac{2}{p}}}{\pi}\max\lbrace\abs{\log(\frac{\rho_1}{2})},\abs{\log(2\sqrt{2}\tilde{\rho}_1)}\rbrace\left(\frac{4m}{\rho_0}+\frac{4\sqrt{\tilde{\rho}_1m}}{\rho_0}+\frac{6\pi C}{\tilde{\rho}_1}\right).
\end{equation}

The final piece of the puzzle is the estimate of
\begin{equation}\label{double integral of Laplacian and Gamma}
\int_{\Omega^{\rho_1}_{\rho_0}}\left(\int_{\Omega^{\tilde{\rho}_1}_{\rho_0}}\abs{ H_N(x,y)\Delta(\alpha-2u)(y)}dy\right)^pdx.
\end{equation}
Here we may use Minkowski's inequality for integrals once more to see that the above is bounded by
\begin{equation}
\left(\int_{\Omega^{\tilde{\rho}_1}_{\rho_0}}\abs{\Delta(\alpha-2u)(y)}\left(\int_{\Omega^{\rho_1}_{\rho_0}}\abs{ H_N(x,y)}^pdx\right)^\frac{1}{p}dy\right)^p.
\end{equation}
Thus, we may bound (\ref{double integral of Laplacian and Gamma}) by
\begin{equation}
	\left(\int_{\Omega^{\tilde{\rho}_1}_{\rho_0}}\abs{\Delta(\alpha-2u)(y)}\left(\int_{\Omega^{\tilde{\rho}_1}_{\rho_0}\cup F\left(\Omega^{\tilde{\rho}_1}_{\rho_0}\right)}\abs{ H_N(x,y)}^pdx\right)^\frac{1}{p}dy\right)^p.
\end{equation}
Again, using the change of variable formula and (\ref{log estimate}), we bound (\ref{double integral of Laplacian and Gamma}) by
\begin{equation}\label{Bound on the Interior}
\left(\left[\frac{1}{2}+16\tilde{\rho}_1^2\log(2\sqrt{2}\tilde{\rho}_1))\right]^{\frac{2-p}{p}}\left[\frac{1}{2\pi}\sqrt{2\pi+64\pi\tilde{\rho}_1^2\log(2\sqrt{2}\tilde{\rho}_1)^2}\right]^{\frac{2p-2}{p}}\frac{4m+4\sqrt{\tilde{\rho}_1m}}{\rho_0}\right)^p.
\end{equation}

Putting everything above together shows that
\begin{equation}
	\int_{\Omega^{\rho_1}_{\rho_0}}\abs{\alpha-2u}^p\le C(p)^2(\ref{total asym flat part of boundary estimate})+C(p)^3\left((\ref{Bound on mu bottom})^p+(\ref{Bound on mu top})^p\right)+C(p)(\ref{Bound on the Interior}).
\end{equation}
Thus, for any $\epsilon>0$ and $\rho_1>\rho_0$ we can pick an appropriate $\tilde{\rho}_1$ and ADM mass $m$ so that
\begin{equation}
\abs{\abs{\alpha-2u}}_{L^1\left(\Omega^{\rho_1}_{\rho_0}\right)}<\frac{\epsilon}{2}.
\end{equation}

We can get similar estimates for $\nabla(\alpha-2u)$ by differentiating the representation formula:
\begin{equation}
\nabla(\alpha-2u)(x)=\int_{\partial\Omega^{\tilde{\rho}_1}_{\rho_0}}(\alpha-2u)\nabla_x\frac{\partial H_N}{\partial\nu}-\nabla_x H_N(x,y)\frac{\partial(\alpha-2u)}{\partial\nu}+\int_{\Omega_{\rho_0}^{\tilde{\rho}_1}}(\nabla_x H_N)\Delta(\alpha-2u).
\end{equation}
We see that
\begin{equation}
\int_{\Omega^{\rho_1}_{\rho_0}}\abs{\nabla(\alpha-2u)}^p\le C(p)\int_{\Omega^{\rho_1}_{\rho_0}}\left(\int_{\partial\Omega^{\tilde{\rho}_1}_{\rho_0}}\abs{\frac{\partial(\alpha-2u)}{\partial\nu}\nabla_xH_N}+\abs{(\alpha-2u)\nabla_x\frac{\partial H_N}{\partial\nu}}\right)^p+\left(\int_{\Omega^{\tilde{\rho}_1}_{\rho_0}}\abs{\nabla_xH_N}\abs{\Delta(\alpha-2u)}\right)^p.
\end{equation}

As before, we will break $\partial\Omega^{\tilde{\rho}_1}_{\rho_0}$ into (\ref{Asymptotically flat part of the boundary}) and (\ref{Interior part of the boundary}). We start with (\ref{Asymptotically flat part of the boundary}). A quick calculation shows that
\begin{equation}
\abs{\nabla_xH_N}\le\frac{1}{2\pi}\left(\frac{1}{\abs{x-y}}+\frac{1}{\abs{\bar{x}-y}}\right)
\end{equation}
and
\begin{equation}
	\abs{\nabla_x\frac{\partial H_N}{\partial\nu}}\le\frac{3}{2\pi}\left(\frac{1}{\abs{x-y}^2}+\frac{1}{\abs{\bar{x}-y}^2}\right).
\end{equation}\
Estimating the integral over (\ref{Asymptotically flat part of the boundary}) now proceeds as before.

As a first step in estimating the integral over (\ref{Interior part of the boundary}), we note that
\begin{equation}
	\nabla_x\frac{\partial H_N}{\partial\nu}\Big|_{\lbrace\rho=\rho_0\rbrace}=0.
\end{equation}
Next, we again break (\ref{Interior part of the boundary}) into (\ref{mu bottom}) and (\ref{mu top}). For both pieces we proceed much as we did before. On (\ref{mu bottom}) it is crucial that $p<2$, since it is only then that the integral
\begin{equation}
	\int_{\Omega^{\rho_1}_{\rho_0}}\abs{\nabla_xH_N}^p
\end{equation}
bounded for all $y$ in (\ref{mu bottom}). For (\ref{mu top}), the necessary changes in the argument are straightforward.

Finally, to estimate
\begin{equation}
\int_{\Omega_{\rho_0}^{\rho_1}}\left(\int_{\Omega^{\tilde{\rho}_1}_{\rho_0}}\abs{(\nabla_x H_N)\Delta(\alpha-2u)}\right)^p\le\int_{\Omega_{\rho_0}^{\tilde{\rho}_1}}\left(\int_{\Omega^{\tilde{\rho}_1}_{\rho_0}}\abs{(\nabla_x H_N)\Delta(\alpha-2u)}\right)^p
\end{equation}
we may use the Riesz potential estimates (\ref{Riesz Potential Estimate}) with the appropriate choice of constants. Thus, for $\tilde{\rho}_1$ chosen large enough and $m$ chosen small enough, we may conclude that
\begin{equation}
\abs{\abs{\alpha-2u}}_{W^{1,p}\left(\Omega^{\rho_1}_{\rho_0}\right)}<\epsilon.\qed
\end{equation}

In the course of proving Proposition (\ref{Sobolev Convergence of 2u-alpha}) we actually proved a little more. For future convenience, we record this result as the following corollary.
\begin{cor}\label{2u-alpha Blowup  Estimate}
	Let $\mathcal{M}$ be a family of axisymmetric metrics with nonnegative scalar curvature which is uniformly asymptotically flat outside of radius $R_0$. Suppose also that $\mathcal{M}$ is radially monotone at $\rho_0$. 
	For any $\rho_1>\rho_0$, $\epsilon>0$, and $1\le p<2$ there exists a $\delta>0$ such that if $g\in\mathcal{M}$ and
	\begin{equation}
	m(g)<\delta,
	\end{equation}
	then
	\begin{equation}
	\int_{\Omega^{\rho_1}_{\rho_0}}\abs{\alpha-2u}^p\le\frac{\epsilon}{\rho_0^p}
	\end{equation}
	and
	\begin{equation}
	\int_{\Omega^{\rho_1}_{\rho_0}}\abs{\nabla(\alpha-2u)}^p\le\frac{\epsilon}{\rho_0^p}.
	\end{equation}
\end{cor}

Having successfully estimated the $W^{1,p}$ norm of $\alpha-2u$, we must now turn to estimating the $W^{1,p}$ norm of $e^{\alpha-2u}$. As was noted earlier, control over the $W^{1,p}$ norm of $\alpha-2u$ for $1\le p<2$ falls short of what we need to apply the Moser-Trudinger inequality to $\alpha-2u$. It is thus not immediately clear how to turn estimates for $\alpha-2u$ into estimates for $e^{\alpha-2u}$. Luckily, the special nature of the fundamental solution to the Laplacian in two dimensions allows us to prove a Moser-Trudinger like inequality which we can use on $\alpha-2u$.
\begin{lemma}\label{Moser Trudinger like Inequality}
	Let $\Omega$ be a  bounded domain in the plane on which the divergence theorem holds and let $\Gamma$ be the fundamental solution for the Laplacian. Suppose we have $\psi\in C^2(\Omega)\cap C^{1}(\bar{\Omega})$ and $\Delta\psi\in L^1(\Omega)$. Let $\Omega_\sigma$ denote $=\lbrace x\in\Omega:d(x,\partial\Omega)\ge\sigma\rbrace$ and let $r_0=\max\lbrace1,diam(\Omega)\rbrace$. Then we have the estimate:
	\begin{equation}
	\int_{\Omega_\sigma}e^{|\psi|}\le\left(|\Omega_\sigma|+\frac{\pi\abs{\abs{\Delta\psi}}_1}{4\pi-\abs{\abs{\Delta\psi}}_1}+2\pi(r_0-1)r_0[r_0^{\frac{\abs{\abs{\Delta\psi}}}{2\pi}}-1]\right)\sup_{x\in\Omega_\sigma}\exp\left(\int_{\partial\Omega}\left|\psi(y)\frac{\partial\Gamma}{\partial\nu}(x,y)\right|+\left|\Gamma(x,y)\frac{\partial\psi}{\partial\nu}(y)\right|dy\right)
	\end{equation}
\end{lemma}
\textbf{Proof}: From Green's representation we have
\begin{equation}\label{use of Green's rep}
\psi(x)=\int_{\partial\Omega}\psi(y)\frac{\partial\Gamma}{\partial\nu}(x,y)-\Gamma(x,y)\frac{\partial\psi}{\partial\nu}(y)dy+\int_{\Omega}\Gamma(x,y)\Delta\psi(y)dy
\end{equation}
Using the representation formula to rewrite $\int_{\Omega_\epsilon}e^{\abs{\psi}}$, we obtain 
\begin{equation}\label{exp(2u-alpha) terms}
\int_{\Omega_\sigma}e^{|\psi(x)|}dx\le\int_{\Omega_\sigma}\exp\left[\int_{\partial\Omega}\left|\psi(y)\frac{\partial\Gamma}{\partial\nu}(x,y)-\Gamma(x,y)\frac{\partial\psi}{\partial\nu}(y)\right|dy\right]\exp\left[\int_{\Omega}\left|\Gamma(x,y)\Delta\psi(y)\right|dy\right]dx
\end{equation}
We bound the first term on the right pointwise by its supremum over $\Omega_\sigma$. Then we may take it outside of the integrand. 
\begin{equation}\label{Green's rep ineq}
\int_{\Omega_\sigma}e^{|\psi(x)|}dx\le\sup_{x\in\Omega_{\sigma}}\exp\left[\int_{\partial\Omega}\left|\psi(y)\frac{\partial\Gamma}{\partial\nu}(x,y)\right|+\left|\Gamma(x,y)\frac{\partial\psi}{\partial\nu}(y)\right|dy\right]\int_{\Omega_{\sigma}}\exp\left[\int_{\Omega}\left|\Gamma(x,y)\Delta\psi(y)\right|dy\right]dx
\end{equation}
We may now concentrate on estimating
\begin{equation}
\int_{\Omega_\sigma}\exp\left[\int_{\Omega}\abs{\Gamma(x,y)\Delta\psi(y)}dy\right]
\end{equation}
The strategy is to expand the above integral using the Taylor series for the exponential function and then bound each term appearing in the expansion:
\begin{equation}
\int_{\Omega_\sigma}\left(e^{\int_{\Omega}\abs{\Gamma(x,y)\Delta(\alpha-2u)(y)}dy}\right)dx=\sum_{k=0}^{\infty}\int_{\Omega_\sigma}\frac{\left(\int_{\Omega}\abs{\Gamma(x,y)\Delta\psi(y)}dy\right)^k}{k!}dx.
\end{equation}
First, recall that the fundamental solution of the Laplacian in two dimensions is given by
\begin{equation}
\frac{1}{2\pi}\log\abs{x-y}
\end{equation} 
Second, after observing that $\Omega_\sigma\subset\Omega$, and pulling constants out, we get the inequality
\begin{equation}\label{Taylor Series of Integral}
\int_{\Omega_\sigma}\frac{\abs{\int_{\Omega}\Gamma(x-y)\Delta\psi(y)dy}^k}{k!}\le\frac{1}{k!(2\pi)^k}\int_{\Omega}\left(\int_{\Omega}\abs{\Delta\psi(y)}\abs{\log(\abs{x-y})}dy\right)^kdx
\end{equation}
We apply Jensen's inequality to the integral on the right to obtain 
\begin{equation}\label{Jensen Inequality Estimate}
\frac{1}{(2\pi)^kk!}\int_{\Omega}\left(\int_{\Omega}\abs{\log(\abs{x-y})}\abs{\Delta\psi(y)}dy\right)^kdx\le\frac{||\Delta\psi||_1^{k-1}}{(2\pi)^kk!}\int_{\Omega}\int_\Omega\abs{\log(\abs{x-y})}^k\abs{\Delta\psi(y)}dydx
\end{equation}
We now use Tonelli's theorem to switch the order of integration to get
\begin{equation}\label{Swap Order of Integration}
\int_{\Omega}\int_\Omega\abs{\log(\abs{x-y})}^k\abs{\Delta\psi(y)}dydx=\int_{\Omega}\abs{\Delta\psi(y)}\int_{\Omega}\abs{\log(\abs{x-y})}^kdxdy
\end{equation}

Putting (\ref{log estimate}), (\ref{Swap Order of Integration}), and (\ref{Jensen Inequality Estimate}) together gives
\begin{equation}
\frac{1}{k!}\int_{\Omega}\abs{\int_{\Omega}\frac{1}{2\pi}\log(\abs{x-y})\Delta\psi(y)dy}^kdx\le\frac{\abs{\abs{\Delta\psi}}_1^k}{(2\pi)^kk!}\left(\frac{\pi k!}{2^k}+2\pi(r_0-1)r_0\log(r_0)^k\right)
\end{equation}
After a quick application of the Monotone Convergence theorem to the summation over $k$ from $k=1$ to infinity of (\ref{Taylor Series of Integral}) we get
\begin{equation}
\int_{\Omega_\sigma}e^{\abs{\int_{\Omega}\Gamma(x,y)\Delta\psi(y)dy}}dx\le|\Omega_\sigma|+\frac{\pi\abs{\abs{\Delta\psi}}_1}{4\pi-\abs{\abs{\Delta\psi}}_1}+(r_0-1)r_0\left[\exp\left(\frac{\log(r_0)\abs{\abs{\Delta\psi}}_1}{2\pi}\right)-1\right].\qed
\end{equation}

We have the following corollary, which is the actual inequality we will use.
\begin{cor}\label{Msr-like corollary}
	Suppose $\psi\in C^2\left(\Omega^{\rho_1}_{\rho_0}\right)\cap C^1\left(cl\left(\Omega^{\rho_1}_{\rho_0}\right)\right)$ and let $r_0=\max\lbrace1,diam\left(\Omega^{\rho_1}_{\rho_0}\right)\rbrace$. Then
	\begin{equation}
	\int_{\left(\Omega^{\rho_1}_{\rho_0}\right)_{\sigma}}e^{\abs{\psi}}
	\end{equation}
	is bounded above by
	\begin{equation}
	e^{C(\sigma,\rho_1)\abs{\abs{\Delta\psi}}_1}\left(\left|\left(\Omega^{\rho_1}_{\rho_0}\right)_\sigma\right|+\frac{\pi\abs{\abs{\Delta\psi}}_1}{4\pi-\abs{\abs{\Delta\psi}}_1}+r^2_0[r_0^{\frac{\abs{\abs{\Delta\psi}}}{2\pi}}-1]\right)\sup_{x\in\left(\Omega^{\rho_1}_{\rho_0}\right)_\sigma}\exp\left(\int_{\partial\Omega}|\psi(y)\frac{\partial H_N}{\partial\nu}|+|H_N\frac{\partial\psi}{\partial\nu}(y)|dy\right),
	\end{equation}
	where $C(\sigma,\rho_1)=\frac{1}{2\pi}\max\lbrace\abs{\log(\sigma)},\abs{\log(2\sqrt{2}\rho_1)}\rbrace$.
\end{cor}
\textbf{Proof:} If we replace $\Gamma$ by $H_N$ in (\ref{use of Green's rep}), then the right hand side of (\ref{Green's rep ineq}) becomes
\begin{equation}
	\sup_{x\in\left(\Omega^{\rho_1}_{\rho_0}\right)_{\sigma}}\exp\left[\int_{\partial\Omega^{\rho_1}_{\rho_0}}\left|\psi(y)\frac{\partial H_N}{\partial\nu}\right|+\left| H_N\frac{\partial\psi}{\partial\nu}(y)\right|dy+\int_{\Omega^{\rho_1}_{\rho_0}}\abs{\Gamma(\bar{x},y)\Delta\psi}\right]\int_{\left(\Omega^{\rho_1}_{\rho_0}\right)_{\sigma}}\exp\left[\int_{\Omega}|\Gamma(x,y)\Delta\psi(y)|dy\right].
\end{equation}
We see that 
\begin{equation}
	\sup_{x\in\left(\Omega^{\rho_1}_{\rho_0}\right)_{\sigma}}\int_{\Omega^{\rho_1}_{\rho_0}}\abs{\Gamma(\bar{x},y)\Delta\psi}\le C(\sigma,\rho_1)\abs{\abs{\Delta\psi}}_{L^1\left(\Omega^{\rho_1}_{\rho_0}\right)}.
\end{equation}
The corollary now follows from Lemma \ref{Moser Trudinger like Inequality}.$\qed$

In order to apply Corollary \ref{Msr-like corollary} to $\alpha-2u$, we need an $L^1$ bound on $\Delta(\alpha-2u)$ and an uniform bound on the boundary. In Lemma \ref{L1 Estimate for the Laplacian of 2u-alpha} we established the necessary $L^1$ bound. Now, we will demonstrate the needed uniform control on the boundary. The following result is very similar to Lemma \ref{Inner Bdry Est for alpha-2u}, however, due to technical necessities, the statement and proof are slightly different.
\begin{lemma}\label{Boundary Estimate for alpha-2u}
	Let $\mathcal{M}$ be a family of axisymmetric metrics with nonnegative scalar curvature which is uniformly asymptotically flat outside of radius $R_0$. Suppose also that $M$ is radially monotone at $\rho_0$. Let $\Omega^{\rho_1}_{\rho_0}$ denote the region
	\begin{equation}
		\lbrace(\rho,z)|\rho_0\le\rho\le \rho_1,\abs{z}\le\frac{\rho_1}{2}\rbrace,
	\end{equation}
	and $(\Omega^{\rho_1}_{\rho_0})_{\sigma}$ denote $\lbrace x\in\Omega^{\rho_1}_{\rho_0}|d(x,\partial\Omega^{\rho_1}_{\rho_0})>\sigma\rbrace$. Let $\rho_1\ge R_0$. If $g\in\mathcal{M}$ and the ADM mass of $g$ is less than $m$,
	then
	\begin{equation}\label{boundary size}
	\sup_{x\in(\Omega^{\rho_1}_{\rho_0})_\sigma}\exp\left(\int_{\partial\Omega^{\rho_1}_{\rho_0}}\left|H_N(x,y)\frac{\partial(2u-\alpha)}{\partial\nu}(y)\right|+\left|(2u-\alpha)(y)\frac{\partial H_N}{\partial\nu}(x,y)\right|dy\right)\le\exp\left[ C(m,\sigma,\rho_1,\rho_0)\right]
	\end{equation}
	where
	\begin{equation}
		C(m,\sigma,\rho_1,\rho_0)=\max\left\lbrace\abs{\log 2\sqrt{2}\rho_1},\abs{\log\sigma}\right\rbrace\left(\frac{4m+4\sqrt{\rho_1 m}}{\pi\rho_0}+\frac{9C}{\rho_1}\right)+\frac{3C}{\sigma}.
	\end{equation}
\end{lemma}
\textbf{Proof}: As we observed earlier, for three sides of the rectangle $\Omega^{\rho_1}_{\rho_0}$, the necessary estimates to control the left hand side of (\ref{boundary size}) follow from the uniformly asymptotically flat condition.  Let's make this more precise. First, consider those pieces of the rectangle parallel to the $\rho$-axis. Here
\begin{equation}
\frac{\partial}{\partial\nu}=\pm\frac{\partial}{\partial z}.
\end{equation}
From the definitions, we know that
\begin{equation}
\abs{\frac{\partial u}{\partial z}(y)}\le\frac{C}{\abs{y}^2}
\end{equation}
and
\begin{equation}
\abs{\frac{\partial \alpha}{\partial z}(y)}\le\frac{C}{\abs{y}^2}	
\end{equation}
for
\begin{equation}
\abs{y}\ge  R_0.
\end{equation}
We may combine these two inequalities using the triangle inequality to conclude that
\begin{equation}
\abs{\frac{\partial(\alpha-2u)}{\partial z}(y)}\le\frac{3C}{\abs{y}^2}.
\end{equation}
Analogously, we have
\begin{equation}
\abs{\alpha-2u}\le\frac{3C}{\abs{y}}.
\end{equation}
In fact, the same is true on the final edge, so the above estimates are true on all of $\partial\Omega^{\rho_1}_{\rho_0}-\lbrace\rho=\rho_0\rbrace$.

Armed with these estimates, let's take a look at the integral
\begin{equation}
\int_{\partial\Omega^{\rho_1}_{\rho_0}-\lbrace\rho=\rho_0\rbrace} \left|H_N(x,y)\frac{\partial(\alpha-2u)}{\partial\nu}(y)\right|+\left|(\alpha-2u)(y)\frac{\partial H_N}{\partial\nu}(x,y)\right|dy
\end{equation}
Since the point $x$ is at a distance of at least $\sigma$ away from the boundary, we know that
\begin{equation}
\frac{\partial H_N}{\partial \nu}\le\frac{1}{\sigma\pi}
\end{equation}
and
\begin{equation}\label{unfm estimate for HN away from bdry}
| H_N(x,y)|\le \frac{1}{\pi}\max\left\lbrace|\log(2\sqrt{2}\rho_1)|,|\log(\sigma)|\right\rbrace
\end{equation}
To start, we can bound
\begin{equation}
\int_{\partial\Omega^{\rho_1}_{\rho_0}-\lbrace\rho=\rho_0\rbrace}\left|(\alpha-2u)(y)\frac{\partial H_N}{\partial\nu}(x,y)\right|dy
\end{equation}
from above by
\begin{equation}
\int_{\partial\Omega^{\rho_1}_{\rho_0}-\lbrace\rho=\rho_0\rbrace}\frac{3C}{\sigma\pi\abs{y}}dy\le \frac{3C}{\sigma},
\end{equation}
since $\abs{y}\ge\rho_1$ for $y$ in $\partial\Omega^{\rho_1}_{\rho_0}-\lbrace\rho=\rho_0\rbrace$.
We now make a similar estimate for
\begin{equation}
\int_{\partial\Omega^{\rho_1}_{\rho_0}-\lbrace\rho=\rho_0\rbrace} \left|H_N(x,y)\frac{\partial(\alpha-2u)}{\partial\nu}(y)\right|dy.
\end{equation}
As we did before, we may bound this quantity from above by
\begin{equation}
\int_{\partial\Omega^{\rho_1}_{\rho_0}-\lbrace\rho=\rho_0\rbrace}\frac{3C}{\pi\abs{y}^2}\max\lbrace\abs{\log 2\sqrt{2}\rho_1},\abs{\log\sigma}\rbrace dy\le\frac{3C}{\rho_1}\max\left\lbrace\abs{\log 2\sqrt{2}\rho_1},\abs{\log\sigma}\right\rbrace.
\end{equation}
We need to estimate
\begin{equation}
	\int_{(\partial\Omega^{\rho_1}_{\rho_0})\cap\lbrace\rho=\rho_0\rbrace}\abs{H_N(x,y)\frac{\partial(\alpha-2u)}{\partial\nu}}
\end{equation}
for $x\in(\Omega^{\rho_1}_{\rho_0})_{\sigma}$.
Using (\ref{unfm estimate for HN away from bdry}) and (\ref{Inner Bdry Est for alpha-2u})
we get
\begin{equation}
	\int_{(\partial\Omega^{\rho_1}_{\rho_0})\cap\lbrace\rho=\rho_0\rbrace}\abs{H_N(x,y)\frac{\partial(\alpha-2u)}{\partial\nu}}\le \frac{1}{\pi}\max\left\lbrace|\log(2\sqrt{2}\rho_1)|,|\log(\sigma)|\right\rbrace\left(\frac{4m}{\rho_0}+\frac{4\sqrt{\rho_1m}}{\rho_0}+\frac{6\pi C}{\rho_1}\right).
\end{equation}
Putting the estimates together gives
\begin{equation}
\sup_{x\in(\Omega^{\rho_1}_{\rho_0})_\sigma}\exp\left(\int_{\partial\Omega^{\rho_1}_{\rho_0}} \left|H_N(x,y)\frac{\partial(\alpha-2u)}{\partial\nu}(y)\right|+\left|(\alpha-2u)(y)\frac{\partial H_N}{\partial\nu}(x,y)\right|dy\right)\le C(m,\sigma,\rho_1,\rho_0).\qed
\end{equation}

With all of the above estimates in hand, controlling the $W^{1,p}$ norm of $e^{\alpha-2u}$ is relatively straightforward. The technical requirements of Corollary \ref{Msr-like corollary} force us to consider regions $\Omega^{\rho_1}_{\rho_0}(\sigma)$ for positive $\sigma$, see (\ref{orbit rectangle}).
\begin{lemma}\label{Sobolev Convergence of exp(2u-alpha)}
	Let $\mathcal{M}$ be a family of axisymmetric metrics with nonnegative scalar curvature which is uniformly asymptotically flat outside of radius $R_0$. Let $\Omega^{\rho_1}_{\rho_0}$ denote the region $\lbrace(\rho,z)|\rho_0\le\rho\le \rho_1,\abs{z}\le\frac{\rho_1}{2}\rbrace$. Suppose that $\mathcal{M}$ is also radially monotone at $\rho_0$.
	For every $\rho_1>\max\lbrace\rho_0,R_0\rbrace$, $\epsilon>0$, $\sigma>0$, and $1\le p<2$ there exists a $\delta>0$ such that if the ADM mass of $g\in\mathcal{M}$ is less than $\delta$, then
	\begin{equation}
	\abs{\abs{e^{\abs{\alpha-2u}}-1}}_{W^{1,p}\left(\Omega^{\rho_1}_{\rho_0}(\sigma)\right)}<\epsilon.
	\end{equation}
\end{lemma}
\textbf{Proof}: By assumption, $\alpha-2u$ is bounded and has bounded derivatives, although we make no assumption on what these bounds might be. Thus, we have that $e^{\abs{\alpha-2u}}$ is Lipschitz
As in Lemma \ref{Sobolev Convergence for exp(u)}, we get
\begin{equation}
	\int_{\Omega^{\rho_1}_{\rho_0}(\sigma)}\abs{\nabla e^{\abs{\alpha-2u}}-1}^p\le\int_{\Omega^{\rho_1}_{\rho_0}(\sigma)}\abs{\nabla(\alpha-2u)}^pe^{p\abs{\alpha-2u}}.
\end{equation}
Let $r>1$ be such that $rp<2$. Applying H{\"o}lder's inequality to the above gives
\begin{equation}
	\left(\int_{\Omega^{\rho_1}_{\rho_0}(\sigma)}\abs{\nabla(\alpha-2u)}^{rp}\right)^{\frac{1}{r}}\left(\int_{\Omega^{\rho_1}_{\rho_0}(\sigma)}e^{r'p\abs{\alpha-2u}}\right)^{\frac{1}{r'}},
\end{equation}
where $r'$ is the conjugate exponent to $r$. In order to control the left hand side we appeal to Proposition \ref{Sobolev Convergence of 2u-alpha}. In order to bound the right hand side we first note that
\begin{equation}
\Omega^{\rho_1}_{\rho_0}(\sigma)\subset\left(\Omega^{\rho_1+\sigma}_{\rho_0}\right)_{\sigma}.
\end{equation}
Thus
\begin{equation}
	\int_{\Omega^{\rho_1}_{\rho_0}(\sigma)}e^{r'p\abs{\alpha-2u}}\le\int_{\left(\Omega^{\rho_1+\sigma}_{\rho_0}\right)_{\sigma}}e^{r'p\abs{\alpha-2u}}.
\end{equation}
We may apply Corollary \ref{Msr-like corollary} to the function $r'p(\alpha-2u)$ and modify Lemma \ref{Boundary Estimate for alpha-2u} as necessary in order to see that
\begin{equation}
	\int_{\left(\Omega^{\rho_1+\sigma}_{\rho_0}\right)_{\sigma}}e^{r'p\abs{\alpha-2u}}
\end{equation}
is uniformly bounded for all $m$ small enough. Thus, combining the two estimates above shows that
\begin{equation}
	\abs{\abs{\nabla e^{\abs{\alpha-2u}}}}_{L^p\left(\Omega^{\rho_1}_{\rho_0}(\sigma)\right)}<\frac{\epsilon}{2}
\end{equation}
for sufficiently small $m$. Similarly, for $m$ small enough, we can show that
\begin{equation}
	\abs{\abs{e^{\abs{\alpha-2u}}}}_{L^p\left(\Omega^{\rho_1}_{\rho_0}(\sigma)\right)}<\frac{\epsilon}{2}. \qed
\end{equation}
\section{Proofs of the Theorems}

In this section we will apply the lemmas to prove the theorems stated in the introduction. Most of the above lemmas analyzed functions over the rectangles $\Omega^{\rho_1}_{\rho_0}$. Now we move our focus to the cylindrical annuli
\begin{equation}\label{Cylindrical Annulus}
	\tilde{\Omega}^{\rho_1}_{\rho_0}(\sigma)=\Omega^{\rho_1}_{\rho_0}(\sigma)\times[0,2\pi),
\end{equation}
see (\ref{cylindrical region}).
Except for the final theorem, this change of focus doesn't involve any new difficulties.
\subsection{Proof of Theorem~\ref{W1p estimates}:}
We first restate the theorem.
\begin{theorem}
	Let $\mathcal{M}$ be a family of axisymmetric metrics with nonnegative scalar curvature which is uniformly asymptotically flat outside of radius $R_0$. Suppose that $\mathcal{M}$ is radially monotone at $\rho_0$ and that for each metric in $\mathcal{M}$, we have
	\begin{equation}
	A=B=0.
	\end{equation}
	For every $\rho_1>\max\lbrace\rho_0,R_0\rbrace$, $\epsilon>0$, $\sigma>0$, and $1\le p<2$ there exists a $\delta>0$ such that if the ADM mass of $g\in\mathcal{M}$ is less than $\delta$, then
	\begin{equation}
	\abs{\abs{g-\delta_{\mathbb{R}^3}}}_{W^{1,p}\left(\tilde{\Omega}^{\rho_1}_{\rho_0}(\sigma)\right)}<\epsilon,
	\end{equation}
	and
	\begin{equation}
	\abs{\abs{q-\delta_{\mathbb{R}^2}}}_{W^{1,p}\left({\Omega}^{\rho_1}_{\rho_0}(\sigma)\right)}<\epsilon,
	\end{equation}
	where $\delta_{\mathbb{R}^3}$ denotes the Euclidean metric in cylindrical coordinates, $\delta_{\mathbb{R}^2}$ denotes the Euclidean metric in the $(\rho,z)$ plane, and $q$ denotes the orbit metric of $g$ in the $(\rho,z)$ plane. $\tilde{\Omega}^{\rho_1}_{\rho_0}(\sigma)$ denotes the cylinder given in (\ref{cylindrical region}) and $\Omega^{\rho_1}_{\rho_0}(\sigma)$ denotes its orbit space.
\end{theorem}

\textbf{Proof}: Since we have assumed that $A=B=0$, in order to show that $g$ is $W^{1,p}$ close to $\delta_{\mathbb{R}^3}$ for small ADM mass, we need only show  that 
\begin{equation}\label{Killing Field Estimate}
	\abs{\abs{\rho^2 e^{-2u}-\rho^2}}_{W^{1,p}\left(\tilde{\Omega}^{\rho_1}_{\rho_0}(\sigma)\right)}<\epsilon
\end{equation}
and
\begin{equation}\label{Orbit Metric Estimate}
\abs{\abs{e^{2\alpha-2u}-1}}_{W^{1,p}\left(\tilde{\Omega}^{\rho_1}_{\rho_0}(\sigma)\right)}<\epsilon
\end{equation}
if the ADM mass is sufficiently small. For (\ref{Killing Field Estimate}) this follows quickly from Lemma \ref{Sobolev Convergence for exp(u)}. Demonstrating (\ref{Orbit Metric Estimate}) is only a little more difficult.

As before, we see that
\begin{equation}
	\int_{\tilde{\Omega}^{\rho_1}_{\rho_0}}\abs{e^{2(\alpha-u)}-1}\le\int_{\tilde{\Omega}^{\rho_1}_{\rho_0}}\abs{2\alpha-2u}^pe^{2p(\alpha-u)}.
\end{equation}
After applying H{\"o}lder's inequality to the above with some $r>1$ such that $rp<2$ we obtain
\begin{equation}
	\left(\int_{\tilde{\Omega}^{\rho_1}_{\rho_0}}\abs{2(\alpha-u)}^{rp}\right)^{\frac{1}{r}}\left(\int_{\tilde{\Omega}^{\rho_1}_{\rho_0}}e^{2pr'(\alpha-u)}\right)^{\frac{1}{r'}}.
\end{equation}
In order to estimate the above, we first observe that
\begin{equation}
	2(\alpha-u)=2u+2(\alpha-2u).
\end{equation}
We can now estimate the left hand term using the triangle inequality, Corollary \ref{Sobolev Convergence for u}, and Proposition \ref{Sobolev Convergence of 2u-alpha} for the exponent $rp<2$. For the right hand side we have
\begin{equation}
	\int_{\tilde{\Omega}^{\rho_1}_{\rho_0}}e^{2pr'(\alpha-u)}=\int_{\tilde{\Omega}^{\rho_1}_{\rho_0}}e^{2pr'u}e^{2pr'(\alpha-2u)}.
\end{equation}
After applying H{\"o}lder's inequality, we may use Lemma \ref{Sobolev Convergence of exp(2u-alpha)} and Lemma \ref{Sobolev Convergence for exp(u)} applied to $2pr'u$ and $2pr'(\alpha-2u)$, respectively, to bound the $L^p$ norm of $e^{2\alpha-2u}$. In fact, in the same way, for any fixed $q$ we can bound the $L^q$ norm of $e^{2\alpha-2u}$ for all $m$ small enough, depending on $\rho_1$, $\rho_0$, and $q$. For what follows, we pick $q$ large enough, depending on $p$.
If we take the gradient of $e^{2\alpha-2u}$ we get
\begin{equation}
(e^{2\alpha-2u})\nabla(2\alpha-2u)=e^{2\alpha-2u}(\nabla 2u+2\nabla(\alpha-2u)).
\end{equation}

We again use H{\"o}lder's inequality, Lemma \ref{Sobolev Convergence for exp(u)}, Proposition \ref{Sobolev Convergence of 2u-alpha} and Lemma \ref{Sobolev Convergence of exp(2u-alpha)} to control the $L^p$ norm of $\nabla e^{2\alpha-2u}$.$\qed$

\subsection{Proof of Theorem~\ref{Volume Estimates}:}

Let us first restate the theorem:
\begin{theorem}
	Let $\mathcal{M}$ be a family of axisymmetric metrics with nonnegative scalar curvature which is uniformly asymptotically flat outside of radius $R_0$. Suppose also that $\mathcal{M}$ is radially monotone at $\rho_0$.
	For any constants $\epsilon>0$, $\sigma>0$, and $\rho_1>\max\lbrace\rho_0,R_0\rbrace$, 
	there exists a $\delta>0$ such that if $g\in\mathcal{M}$ and
	\begin{equation}
	m(g) < \delta,
	\end{equation} 
	then
	\begin{equation}
	\abs{\Omega}+\epsilon\ge vol_g(\Omega)\ge\abs{\Omega}-\epsilon
	\end{equation}
	for any region $\Omega$ such that
	\begin{equation}
	\Omega \subset \tilde{\Omega}^{\rho_1}_{\rho_0}(\sigma).
	\end{equation}
\end{theorem}

\textbf{Proof}: A quick calculation shows that the volume form of $g$ in cylindrical coordinates is
\begin{equation}
	\rho e^{2\alpha-3u}d\rho dzd\phi.
\end{equation}
Thus, we have that
\begin{equation}
	\abs{vol_g(\Omega)-|\Omega|}=\abs{\int_{\Omega} \left(e^{2\alpha-3u}-1\right)\rho d\rho dz d\phi}\le\int_{\tilde{\Omega}^{\rho_1}_{\rho_0}}\abs{e^{2\alpha-3u}-1}\rho d\rho dz d\phi.
\end{equation}
As we have done before, we can see that
\begin{equation}
	\int_{\tilde{\Omega}^{\rho_1}_{\rho_0}}\abs{e^{2\alpha-3u}-1}\rho d\rho dz d\phi\le\int_{\tilde{\Omega}^{\rho_1}_{\rho_0}}\abs{2\alpha-3u}e^{\abs{2\alpha-3u}}\rho d\rho dz d\phi.
\end{equation}
We may now apply H{\"o}lder's inequality to the above in order to see that
\begin{equation}
	\int_{\tilde{\Omega}^{\rho_1}_{\rho_0}}\abs{2\alpha-3u}e^{\abs{2\alpha-3u}}\le\left(\int_{\tilde{\Omega}^{\rho_1}_{\rho_0}}\abs{2\alpha-3u}^p\right)^{\frac{1}{p}}\left(\int_{\tilde{\Omega}^{\rho_1}_{\rho_0}}e^{p'\abs{2\alpha-3u}}\right)^{\frac{1}{p'}},
\end{equation}
where $p$ and $p'$ are conjugate exponents and $1\le p<2$. We may use the triangle inequality to make the estimate
\begin{equation}
	\left(\int_{\tilde{\Omega}^{\rho_1}_{\rho_0}}\abs{2\alpha-3u}^p\right)^{\frac{1}{p}}\le \abs{\abs{u}}_{W^{1,p}}+2\abs{\abs{\alpha-2u}}_{W^{1,p}}.
\end{equation}
We may combine Corollary \ref{Sobolev Convergence for u} and Proposition \ref{Sobolev Convergence of 2u-alpha} to control the above. For the exponential term, we use the estimate
\begin{equation}
	e^{p'\abs{2\alpha-3u}}\le e^{p'\abs{u}}e^{2p'\abs{\alpha-2u}}
\end{equation}
and H{\"o}lder's inequality once more to see that
\begin{equation}
	\int_{\tilde{\Omega}^{\rho_1}_{\rho_0}}e^{p'\abs{2\alpha-3u}}\le\left(e^{2p'\abs{u}}\right)^{\frac{1}{2}}\left(\int_{\tilde{\Omega}^{\rho_1}_{\rho_0}}e^{4p'\abs{\alpha-2u}}\right)^{\frac{1}{2}}.
\end{equation}
We now wish to apply Lemma \ref{Sobolev Convergence for exp(u)} and \ref{Sobolev Convergence of exp(2u-alpha)} to the above to see that it is uniformly bounded for $m$ small enough, depending on $\rho_1$, $\rho_0$ and $p$. Combining the two estimates finishes the proof. $\qed$

\subsection{Proof of Theorem~\ref{Area Estimates}:}

Let us first restate the theorem.
\begin{theorem}
	Let $\mathcal{M}$ be a family of axisymmetric metrics with nonnegative scalar curvature which is  uniformly asymptotically flat outside of radius $R_0$. Suppose also that $\mathcal{M}$ is radially monotone at $\rho_0$.
	For any fixed axisymmetric surface $\Sigma$, constant $\epsilon>0$, and constant $\rho_1>\max\lbrace\rho_0,R_0\rbrace$, there exists a $\delta>0$ such that if $m(g)<\delta$, then
	\begin{equation}
	\abs{\Sigma\cap\tilde{\Omega}^{\rho_1}_{\rho_0}(\sigma)}+\epsilon
	\ge Area_g\left(\Sigma\cap\tilde{\Omega}^{\rho_1}_{\rho_0}(\sigma)\right)
	\ge\abs{\Sigma\cap\tilde{\Omega}^{\rho_1}_{\rho_0}(\sigma)}-\epsilon.
	\end{equation}
\end{theorem}

\textbf{Proof}: Let $s$ be a fixed curve in the $(\rho,z)$ plane representing an axisymmetric surface, which we will call $\Sigma$. A calculation shows that the area form associated with $\Sigma$ is
\begin{equation}\label{metric area form}
	\rho\circ s(t)e^{(\alpha-2u)\circ s}\abs{\dot{s}}_{\delta}dtd\phi.
\end{equation}
Note that the Euclidean area form for $\Sigma$ is
\begin{equation}\label{Euclidean area form}
	\rho\circ s(t)\abs{\dot{s}}_{\delta}dtd\phi.
\end{equation} 
From Lemma \ref{Sobolev Convergence of exp(2u-alpha)} we deduce that for any $\epsilon>0$
\begin{equation}
	\abs{\abs{\rho e^{\alpha-2u}-\rho}}_{W^{1,1}\left(\Omega^{\rho_1}_{\rho_0}(\sigma)\right)}<\epsilon,
\end{equation}
if the ADM mass is small enough. Now, the curve segment $s\cap\Omega^{\rho_1}_{\rho_0}(\sigma)$ is part of the boundary of some region in $\Omega^{\rho_1}_{\rho_0}(\sigma)$. Thus, we may use the trace inequality \cite{Evans Gariepy} to conclude that
\begin{equation}
	\abs{\abs{\rho e^{\alpha-2u}-\rho}}_{L^1\left(s\cap\Omega^{\rho_1}_{\rho_0}\right)}<\epsilon.
\end{equation}
This proves the theorem.$\qed$

If the family $\mathcal{M}$ is area enlarging everywhere, then we also have a stronger lower bound on the area of axisymmetric surfaces than the one given above.
\begin{prop}\label{Area Enlarging}
	Let $g$ be an axisymmetric metric. Let $(\rho,z,\phi)$ be the cylindrical coordinates for $g$, let $\delta_{\mathbb{R}^3}$ be the flat metric in cylindrical coordinates, and let $\Sigma$ be a $C^{1}$ axisymmetric surface. If $g$ is area enlarging, then we have
	\begin{equation}\label{area-enl ineq}
		Area_g(\Sigma)\ge Area_{\delta_{\mathbb{R}^3}}(\Sigma)
	\end{equation}
\end{prop}
\textbf{Proof}: Let $\Sigma$ be a $C^1$ axisymmetric surface. Let $s(t)$ be the $C^1$ curve in the $(\rho,z)$ plane which, when revolved around the $\rho$-axis, gives $\Sigma$. We get the following map
\begin{equation}
(t,\phi)\rightarrow(s(t),\phi)
\end{equation}
from $I\times[0,2\pi)$ to $\Sigma$. Let $A_g$ denote the area form of the surface with respect to the metric induced by $g$, and let $A_{\delta_{\mathbb{R}^3}}$ denote the area form induced by the background Euclidean metric. Then using (\ref{metric area form}) and (\ref{Euclidean area form}) we see that 
\begin{equation}
A_{g}=e^{\alpha-2u}A_{\delta_{\mathbb{R}}^3}.
\end{equation}
In coordinates, the area enlarging condition is equivalent to the nonnegativity of $\alpha-2u$. Thus, we know that $e^{\alpha-2u}$ is greater than 1. The result now follows. $\qed$

We may combine the well known Penrose Inequality with the above proposition to constrain the location of outer most minimal surfaces.
\begin{cor}\label{loc outermost min surface}
	Let $\mathcal{M}$ be a family of uniformly asymptotically flat metrics. Suppose $\mathcal{M}$ is either radially monotone or area enlarging. Let $g$ be a metric in $\mathcal{M}$ and $\Sigma$ be the outermost minimal surface. If $\Sigma$ is axisymmetric and topologically a sphere, and
	\begin{equation}
		m(g)\le m,
	\end{equation}
	then
	\begin{equation}
		\Sigma\subset\rho^{-1}\left([0,2\sqrt{2}m)\right).
	\end{equation}
\end{cor}
\textbf{Proof}: Let
\begin{equation}
\rho_0=\max\lbrace\rho:(\rho,z)\in\Sigma\rbrace,
\end{equation}
let $x_0$ be a point in $\Sigma$ point at which $\rho$ attains the maximum $\rho_0$, and let $[x_0]$ denote its orbit under the killing field.
From the Penrose Inequality, we know that
\begin{equation}
m\ge\sqrt{\frac{Area_g(\Sigma)}{16\pi}}.
\end{equation}
Since $\Sigma$ is axisymmetric and topologically a sphere, it must be represented in the $(\rho,z)$ plane by a curve $\gamma$ which intersects the axis of symmetry twice. In particular, $\gamma$ must emanate from the axis, then touch the point $[x_0]$ and then make its way back to the axis. Let $D_{x_0}$ denote the disk represented by a line connecting the axis to the point $[x_0]$. Since this disk has minimal Euclidean area among axisymmetric surfaces with boundary $[x_0]$, we may conclude that
\begin{equation}
	Area_{\delta_{\mathbb{R}^3}}(\Sigma)>2Area_{\delta_{\mathbb{R}^3}}\left(D_{x_0}\right)=2\pi\rho^2_0.
\end{equation} 
Thus, combining the Penrose inequality with the above and the area enlarging inequality (\ref{area-enl ineq}) gives
\begin{equation}
	m>\frac{\rho_0}{2\sqrt{2}}. \qed
\end{equation}

If the metric $g$ in the above has positive scalar curvature, then it is a well known result that the outermost minimal surface must be a sphere. The author does not know if in an axisymmetric metric an outermost minimal surface must also be axisymmetric, though it does seem plausible.

\subsection{Proof of Theorems~\ref{d estimates} and \ref{d estimates pointwise}}
\begin{theorem}
	Let $\mathcal{M}$ be a family of axisymmetric metrics with nonnegative scalar curvature which is uniformly asymptotically flat outside of radius $R_0$. Suppose $\mathcal{M}$ is also radially monotone at $\rho_0$. Additionally, assume that $A=B=0$ in the coordinate representations of the metrics under consideration. Suppose we are given $\epsilon>0$, $\sigma>0$, and $\rho_1>\max\lbrace\rho_0,R_0\rbrace$. There exists a constant $\delta>0$ such that if $m(g)\le\delta$ and $x$ and $y$ are any points such that the Euclidean line segment connecting them lies in $\Omega^{\rho_1}_{\rho_0}(\sigma)\times\lbrace\phi_0\rbrace$ for any $\phi_0$, then 
	\begin{equation}
	d_g(x,y)\le d(x,y)+\epsilon
	\end{equation}
\end{theorem}
\textbf{Proof}: We use the extension theorem for Sobolev functions, appearing as Theorem 4.7 in \cite{Evans Gariepy}. Following the notation of \cite{Evans Gariepy}, we let $U=\Omega^{\rho_1}_{\rho_0}(\sigma)$, $V=2\Omega^{\rho_1}_{\rho_0}(\sigma)$, and $p=1$. Here $V=2\Omega^{\rho_1}_{\rho_0}(\sigma)$ is centered about $U=\Omega^{\rho_1}_{\rho_0}(\sigma)$. Using the extension theorem \cite{Evans Gariepy}, we may see that there is a constant $K$, depending on $\Omega^{\rho_1}_{\rho_0}(\sigma)$, and extensions of the functions $e^{\alpha-2u}-1$, also denoted $e^{\alpha-2u}-1$, such that
\begin{equation}
	\abs{\abs{e^{\alpha-u}-1}}_{W^{1,1}(\mathbb{R}^2)}\le K\abs{\abs{e^{\alpha-u}-1}}_{W^{1,1}\left(\Omega^{\rho_1}_{\rho_0}(\sigma)\right)}.
\end{equation}

In order to obtain an upper estimate for $d_g(x,y)$, it suffices to estimate the length of one curve connecting the points $x$ and $y$. Let $\gamma_{xy}$ denote the Euclidean line in $\Omega^{\rho_1}_{\rho_0}(\sigma)\times\lbrace\phi_0\rbrace$ connecting $x$ to $y$ parameterized by Euclidean arc length In orbit space
\begin{equation}
	\gamma_{xy}(t)=\left(\gamma^{\rho}_{xy}(t),\gamma^{z}_{xy}(t)\right).
\end{equation}
Every such curve lies on the boundary of a square of side length the diameter of $\Omega^{\rho_1}_{\rho_0}(\sigma)$. All such squares are rotations or translations of each other. Thus, there exists a single constant $C$ such that if $\Omega$ is a square with side length the diameter of $\Omega^{\rho_1}_{\rho_0}(\sigma)$, then the trace inequality holds with constant $C$:
\begin{equation}
	\abs{\abs{\omega}}_{L^1(\partial\Omega)}\le C\abs{\abs{\omega}}_{W^{1,1}(\Omega)}.
\end{equation}
Let $l_g(\gamma)$ be the length of $\gamma$ as measured in the metric $g$. Then we have
\begin{equation}
	l_g(\gamma)=\int^{d(x,y)}_{0}e^{(\alpha-u)\circ\gamma(t)}dt.
\end{equation}
We now use the trace inequality \cite{Evans Gariepy} to see that
\begin{equation}
	\abs{d(x,y)-l_g(\gamma)}\le\int^{d(x,y)}_{0}\abs{e^{(\alpha-u)\circ\gamma(t)}-1}dt\le\int_{\partial\Omega}\abs{e^{\alpha-u}-1}\le C\abs{\abs{e^{\alpha-u}-1}}_{W^{1,1}(\Omega)},
\end{equation}
where $\gamma$ lies on the boundary of $\Omega$.
Furthermore, we have
\begin{equation}
	\abs{\abs{e^{\alpha-u}-1}}_{W^{1,1}(\Omega)}\le\abs{\abs{e^{\alpha-u}-1}}_{W^{1,1}(\mathbb{R}^2)}\le K\abs{\abs{e^{\alpha-u}-1}}_{W^{1,1}\left(\Omega^{\rho_1}_{\rho_0}(\sigma)\right)}.
\end{equation}
We may now use Theorem \ref{W1p estimates} to conclude that
\begin{equation}
	\abs{{d(x,y)-l_g(\gamma)}}<\epsilon
\end{equation}
for small enough ADM mass.  $\qed$

Very similarly, we can prove a pointwise upper bound on $d_g(x,y)$ for more general $x$ and $y$ in $\tilde{\Omega}^{\rho_1}_{\rho_0}$.
\begin{theorem}
	Let $\mathcal{M}$ be a family of axisymmetric metrics with nonnegative scalar curvature which is uniformly asymptotically flat outside of radius $R_0$. Suppose also that $\mathcal{M}$ is radially monotone at $\rho_0$. Additionally, assume that $A=B=0$ in the coordinate representations of the metrics under consideration. Suppose we are given $\rho_1>\max\lbrace\rho_0,R_0$, $\epsilon>0$ and $\sigma>0$ and points $x$ and $y$ such that the Euclidean line segment connecting them lies in $\tilde{\Omega}^{\rho_1}_{\rho_0}(\sigma)$. There exists a constant $\delta>0$ such that if $m(g)\le\delta$, then 
	\begin{equation}
	d_g(x,y)\le d(x,y)+\epsilon
	\end{equation}
\end{theorem}
\textbf{Proof:} As before, let $\gamma$ be the Euclidean line connecting $x$ to $y$. Then we have that
\begin{equation}
\abs{l_g(\gamma_{xy})-1}\le\int^{d(x,y)}_{0}\abs{\sqrt{e^{2(\alpha-u)\circ\gamma}\left((\gamma'_{\rho})^2+(\gamma'_z)^2\right)+\gamma_{\rho}^2e^{-2u\circ\gamma}(\gamma'_{\phi})^2}-1}dt.
\end{equation}
Let
\begin{equation}
	Z=e^{\alpha-u}\left(\gamma'_\rho\frac{\partial}{\partial\rho}+\gamma'_{z}\frac{\partial}{\partial z}\right)+e^{-u}\gamma'_{\phi}\frac{\partial}{\partial\phi}.
\end{equation}
Using the reverse triangle inequality, we observe that
\begin{equation}
	\abs{\abs{Z}-1}=\abs{\abs{Z}-\abs{\gamma'}}\le\abs{Z-\gamma'},
\end{equation}
where we are working with the Euclidean metric in cylindrical coordinates. Thus,
we may estimate the above integral by
\begin{equation}
\int^{d(x,y)}_{0}\sqrt{\left(e^{(\alpha-u)\circ\gamma}-1\right)^2\left((\gamma'_{\rho})^2+(\gamma'_z)^2\right)+\left(e^{-u\circ\gamma}-1\right)^2\gamma^2_{\rho}(\gamma'_\phi)^2}dt.
\end{equation}
Using the triangle inequality and the bounds
\begin{equation}
(\tilde{\gamma}'_{\rho})^2+(\tilde{\gamma}'_z)^2\le 1,
\end{equation}
and
\begin{equation}
	\abs{\gamma_{\rho}\gamma'_{\phi}}\le 1,
\end{equation}
we see that the above is bounded in turn by
\begin{equation}
\int^{d(x,y)}_{0}\abs{e^{(\alpha-u)\circ\gamma}-1}dt+\int_{0}^{d(x,y)}\abs{e^{-u\circ\gamma}-1}dt.
\end{equation}
Let $\tilde{\gamma}$ be the projection of $\gamma$ to the $(\rho,z)$ plane. $\tilde{\gamma}$ lies in the boundary of a region $\Omega$. Since $u$ and $\alpha$ don't depend on $\phi$, we see that $u\circ\gamma=u\circ\tilde{\gamma}$ and $\alpha\circ\gamma=\alpha\circ\tilde{\gamma}$. We can now use the trace theorem, and then apply Theorem \ref{W1p estimates} as we did before to show that for ADM mass small enough, we have
\begin{equation}
\int^{d(x,y)}_{0}\abs{e^{(\alpha-u)\circ\tilde{\gamma}}-1}dt+\int_{0}^{d(x,y)}\abs{e^{-u\circ\tilde{\gamma}}-1}dt<\epsilon.\qed
\end{equation}

\subsection{Proof of Theorem~\ref{ufm con. in the ext.}}
We restate the theorem.
\begin{theorem}
	Let $\mathcal{M}$ be a family of axisymmetric metrics with nonnegative scalar curvature which is uniformly asymptotically flat outside of radius $R_0$. Suppose that $\mathcal{M}$ is radially monotone and that for all $g\in\mathcal{M}$ we have 
	\begin{equation}
	A=B=0.
	\end{equation}
	Let $R_1>R_0$ and let $A(R_0,R_1)$ denote the coordinate spherical annulus centered at the origin. For any given $0<\beta<1$ and $\epsilon>0$ there exists a $\delta>0$ such that if $g\in \mathcal{M}$ and
	\begin{equation}
	m(g)<\delta,
	\end{equation}
	then
	\begin{equation}
	\abs{\abs{g-\delta_{\mathbb{R}^3}}}_{C^{0,\beta}\left(A(R_0,R_1)\right)}<\epsilon.
	\end{equation}
\end{theorem}
\textbf{Proof:} Since we have assumed that $A=B=0$, the proof will be established if we can show that 
\begin{equation}
	\abs{\abs{e^{2\alpha-2u}-1}}_{C^{0,\beta}\left(A(R_0,R_1)\right)}<\epsilon
\end{equation}
and
\begin{equation}
	\abs{\abs{e^{-2u}-1}}_{C^{0,\beta}\left(A(R_0,R_1)\right)}<\epsilon
\end{equation}
for small enough ADM mass. The above inequalities will follow if we can show that
\begin{equation}
	\abs{\abs{\alpha-u}}_{C^{0,\beta}\left(A(R_0,R_1)\right)}<\tilde{\epsilon}
\end{equation}
and
\begin{equation}
	\abs{\abs{u}}_{C^{0,\beta}\left(A(R_),R_1\right)}<\tilde{\epsilon}
\end{equation}
for small enough ADM mass, where $\tilde{\epsilon}$ depends on $\epsilon$ above. Using the triangle inequality, we see that it is sufficient to bound the $C^{0,\beta}$ norms of $u$ and  $\alpha-2u$. These bounds are the content of Lemma \ref{Uniform Control of u} and Lemma \ref{Uniform Control of alpha-2u} below, respectively. $\qed$

\begin{lemma}\label{Uniform Control of u}
	Suppose $\mathcal{M}$ is a collection of axisymmetric metrics which is uniformly asymptotically flat outside a ball of radius $R_0$. Let $u$ be the function appearing in the axisymmetric coordinate representation of $g$. Let $R_1$ be greater than $R_0$ and $A(R_0,R_1)$ be the spherical annulus centered at the origin. For $\epsilon>0$ and $0<\beta_0<1$ there exists a $\delta>0$ such that if $g\in\mathcal{M}$ and
	\begin{equation}
		m(g)<\delta,
	\end{equation}
	then
	\begin{equation}\label{uniform estimate for u}
	\abs{\abs{u}}_{C^{0,\beta}\left(A(R_0,R_1)\right)}<\epsilon.
	\end{equation}
\end{lemma}
\textbf{Proof}: Since we are working in the asymptotically flat regime, we have uniform upper bounds on the $C^{1}(A(R_0,R_1))$ norms of the metric functions. From Lemma \ref{estimate for the Sobolev Norm of u in R3} we may bound the $W^{1,2}(A(R_0,R_1))$ norm of $u$. We now interpolate between these two estimates to bound the $W^{1,q}$ norm of $u$ for arbitrarily large $q$. Specifically, we write
\begin{equation}
\int_{A(R_0,R_1)}u^{q}=\int_{A(R_0,R_1)}u^{2}u^{q-2}\le\abs{\abs{u}}^{q-2}_{\infty}\int_{A(R_0,R_1)}u^2
\end{equation}
We may do the same for the derivatives of $u$. In the end, we get the following bounds
\begin{equation}
\abs{\abs{u}}_q\le\abs{\abs{u}}_2^{\frac{2}{q}}\abs{\abs{u}}_\infty^{1-\frac{2}{q}}
\end{equation}
and
\begin{equation}
\abs{\abs{\nabla u}}_q\le\abs{\abs{\nabla u}}_2^{\frac{2}{q}}\abs{\abs{\nabla u}}_{\infty}^{1-\frac{2}{q}}.
\end{equation}
By assumption $\abs{\abs{u}}_{\infty}+\abs{\abs{\nabla u}}_{\infty}\le C$. Furthermore, by Lemma \ref{estimate for the Sobolev Norm of u in R3}, we know $\abs{\abs{u}}_{W^{1,2}(A(R_0,R_1))}<\tilde{\epsilon}$ for sufficiently small $m$. Thus, we obtain the estimate
\begin{equation}
\abs{\abs{u}}_{W^{1,q}}\le C^{1-\frac{2}{q}}\tilde{\epsilon}^{\frac{2}{q}}.
\end{equation}
We may now choose $q$ large enough and appeal to the Sobolev Embedding Theorem to get $C^{0,\beta_0}$ bounds on $u$ for $\beta_0<1$.
$\qed$
\begin{remark}
	It is important to note that we didn't use the hypothesis of radial monotonicity in the above. We only need radial monotonicity to control $\alpha-2u$.
\end{remark}
We will try to produce similar uniform estimates for $\alpha-2u$. However, as before, the process is harder. Whereas for $u$ we started off with $W^{1,p}_{loc}(\mathbb{R}^3)$ control, for $\alpha-2u$ we only have $W^{1,p}_{loc}(\mathbb{R}^2_{+})$ control. Even worse, the estimates we were able to prove become weaker as we approach the axis $\lbrace\rho=0\rbrace$, see Corollary \ref{2u-alpha Blowup  Estimate}. In order to work our way around this conundrum, we must use the extra factor of $\rho$ present in integrating over $B_R$ in $\mathbb{R}^3$ to control the bad behavior seen in Corollary \ref{2u-alpha Blowup  Estimate}.
\begin{lemma}\label{Countering Blowup Estimate}
	Let $f$ be a measurable function on $\Omega^{\rho_1}_{0}$. Suppose for each $t$ we have the estimate
	\begin{equation}
	\int_{\Omega^{\rho_1}_{t}}\abs{f}\le\frac{\epsilon}{t^{q}}
	\end{equation}
	for some $\epsilon>0$ and $q>0$. Suppose $\sigma>q$. Then, there exists a constant, denoted $C(\sigma,q)$, depending only on $\sigma$ and $q$ such that
	\begin{equation}
	\int_{\Omega^{\rho_1}_{0}}\rho^{\sigma}\abs{f}\le C(\sigma,q)\epsilon.
	\end{equation}
\end{lemma}
\textbf{Proof}: Let $t_n=2^{-n}\rho_1$ and 
let $\Omega_{t_n,t_{n-1}}$ be the following rectangle.
\begin{equation}
\Omega_{t_n,t_{n-1}}=\lbrace t_n<\rho\le t_{n-1},\abs{z}\le\frac{\rho_1}{2}\rbrace
\end{equation}
From the Monotone Convergence Theorem we see that
\begin{equation}
\int_{\Omega^{\rho_1}_{\rho_0}}\rho^\sigma\abs{f}^p=\int_{\Omega_{0,t_0}}\rho^{\sigma}\abs{f}^p=\sum_{1}^{\infty}\int_{\Omega_{t_n,t_{n-1}}}\rho^{\sigma}\abs{f}^p.
\end{equation}
We now make the estimate
\begin{equation}
\int_{\Omega_{t_n,t_{n-1}}}\rho^{\sigma}\abs{f}^p\le t^{\sigma}_{n-1}\frac{\epsilon}{t^{q}_{n}}=2^{\sigma}\rho^{\sigma-q}_1(2^{\sigma-q})^{-n}\epsilon.
\end{equation}
This gives a convergent series so long as $\sigma>q$. In total, we have the estimate
\begin{equation}
\int_{\Omega_{0,t_0}}\rho^{\sigma}\abs{f}^p\le C(\sigma,q)\epsilon.\qed
\end{equation}

We now make use of the above lemma to control the $W^{1,1}$ norm of $\alpha-2u$ over the ball of radius $R$ about the origin in $\mathbb{R}^3$.
\begin{lemma}\label{2u-alpha Sobolev Estimate over B_R}
	Let $\mathcal{M}$ be a family of axisymmetric metrics which is uniformly asymptotically flat outside of radius $R_0$. Suppose that $\mathcal{M}$ is also a radially monotone family of metrics. For any $R$ and $\epsilon>0$ there exists a $\delta>0$ such that if $g\in\mathcal{M}$ and
	\begin{equation}
	m(g)<\delta,
	\end{equation}
	then
	\begin{equation}
	\abs{\abs{\alpha-2u}}_{W^{1,1}(B_R)}<\epsilon.
	\end{equation}
\end{lemma}
\textbf{Proof}: Let $D_R$ be the two dimensional half disk of radius $R$ about the origin. Then 
\begin{equation}
\int_{B_R}\abs{\alpha-2u}=2\pi\int_{D_R}\rho\abs{\alpha-2u}
\end{equation}
and
\begin{equation}
\int_{B_R}\abs{\nabla(\alpha-2u)}=2\pi\int_{D_R}\rho\abs{\nabla(\alpha-2u)}.
\end{equation}
For some $\mu>0$, to be specified later, we rewrite the first quantity as
\begin{equation}
\int_{D_R}\rho^{-\mu}\rho^{1+\mu}\abs{\alpha-2u}.
\end{equation}
Let $1<q<2$ and $q'$ be conjugate exponents. We apply H{\"o}lder's inequality to the above to get
\begin{equation}
\left(\int_{D_R}\rho^{-\mu q'}\right)^{\frac{1}{q'}}\left(\int_{D_R}\rho^{(1+\mu)q}\abs{\alpha-2u}^q\right)^{\frac{1}{q}}.
\end{equation}
Choose $\mu$ small enough that
\begin{equation}
\mu q'<1.
\end{equation}
We may pick large enough that $D_R\subset\Omega^{\rho_1}_0$. From Corollary \ref{2u-alpha Blowup  Estimate} and Lemma \ref{Countering Blowup Estimate}, we see that for some constant $C(\mu,q,R)$,
\begin{equation}
\int_{D_R}\rho\abs{\alpha-2u}^p\le C(\mu,q,R)\epsilon
\end{equation}
if $m$ is chosen small enough.The same argument can be made for
\begin{equation}
\int_{D_R}\rho\abs{\nabla(\alpha-2u)}. \qed
\end{equation}

We now make an estimate on the uniform norm of $\alpha-2u$ similar to Lemma \ref{Uniform Control of u}.
\begin{lemma}\label{Uniform Control of alpha-2u}
	Suppose $\mathcal{M}$ is a collection of axisymmetric metrics which is uniformly asymptotically flat outside a ball of radius $R_0$. Let $R_1$ be greater than $R_0$ and $A(R_0,R_1)$ be the spherical annulus centered at the origin. For $\epsilon>0$ and $0<\beta<1$ there exists a $\delta>0$ such that if $g\in\mathcal{M}$ and
	\begin{equation}
	m(g)<\delta,
	\end{equation}
	then
	\begin{equation}\label{uniform estimate for alpha-2u}
	\abs{\abs{\alpha-2u}}_{C^{0,\beta}\left(A(R_0,R_1)\right)}<\epsilon.
	\end{equation}
\end{lemma}
\textbf{Proof:} We imitate the proof of Lemma \ref{Uniform Control of u}. As before, we write
\begin{equation}
	\int_{A\left(R_0,R_1\right)}\abs{\alpha-2u}^q\le\abs{\abs{\alpha-2u}}^{q-1}_{\infty}\int_{A\left(R_0,R_1\right)}\abs{\alpha-2u}.
\end{equation}
We also have
\begin{equation}
	\int_{A\left(R_0,R_1\right)}\abs{\nabla(\alpha-2u)}^q\le\abs{\abs{\nabla(\alpha-2u)}}^{q-1}_{\infty}\int_{A\left(R_0,R_1\right)}\abs{\nabla(\alpha-2u)}.
\end{equation}
By the asymptotic flatness assumption, we know that
\begin{equation}
	\abs{\abs{\alpha-2u}}_{\infty}+\abs{\abs{(\alpha-2u)}}_{\infty}\le C
\end{equation}
For some $C$ depending only on the uniform falloff in Definition \ref{defn-unif-asym}. Thus, for any exponent $q$ we can use Lemma \ref{2u-alpha Sobolev Estimate over B_R} to control the Sobolev norm $\abs{\abs{\alpha-2u}}_{W^{1,q}\left(A(R_0,R_1)\right)}$ by the ADM mass. Using the Sobolev embedding theorem, we see that
\begin{equation}
	\abs{\abs{\alpha-2u}}_{C^{0,\beta}}\le C\abs{\abs{\alpha-2u}}^{\frac{1}{q}}_{W^{1,1}\left(A(R_0,R_1)\right)},
\end{equation}
where $\beta=1-\frac{3}{q}$, the constant $C$ depends only on the uniform falloff in Definition \ref{defn-unif-asym}, the region $A(R_0,R_1)$, and $q$. Now we can use Lemma \ref{2u-alpha Sobolev Estimate over B_R} to control the uniform norm $\alpha-2u$ on $A\left(R_0,R_1\right)$. $\qed$

\section{Area Enlarging Case}
We now show that all the theorems stated hold when we assume our family of uniformly asymptotically flat metrics is area enlarging and strongly uniformly asymptotically flat, instead of radially monotone. The only steps required are to prove a lemma analogous to \ref{Inner Bdry Est for alpha-2u} and a proposition analogous to \ref{Sobolev Convergence of 2u-alpha}. The main difference between the radially monotone case and the area enlarging one is in the choice of function for Green's representation formula. Instead of working with $H_N(x,y)$, we will use $H_D(x,y)$ (\ref{Dirichlet rep}). We also focus on slightly different rectangles,
\begin{equation}
	\Omega^{L}_{\rho_0\rho_1}:=\lbrace(\rho,z):\rho_0\le\rho\le\rho_1,\abs{z}\le\frac{L}{2}\rbrace.
\end{equation}
We now prove the first key lemma for the area enlarging and strongly uniformly asymptotically flat case.
\begin{lemma}\label{2u-alpha boundary estimate}
	Let $\mathcal{M}$ be a family of axisymmetric metrics which is strongly uniformly asymptotically flat outside of radius $R_0$. Suppose also that $\mathcal{M}$ is area enlarging at $\rho_0$. For any $\rho_1>\rho_0$, $L>0$, and $\epsilon>0$ there exists a $\delta>0$ such that if 
	\begin{equation}
	m(g)<\delta,
	\end{equation}
	then
	\begin{equation}
	\int_{\partial\Omega^{L}_{\rho_0\rho_1}\cap\lbrace\rho=\rho_0\rbrace}\abs{\alpha-2u}<\epsilon.
	\end{equation}
\end{lemma}

Proof: Observe that if $\tilde{L}>L$, then
\begin{equation}
	\int_{\partial\Omega^{\tilde{L}}_{\rho_0\rho_1}\cap\lbrace\rho=\rho_0\rbrace}\abs{\alpha-2u}\ge\int_{\partial\Omega^{L}_{\rho_0\rho_1}\cap\lbrace\rho=\rho_0\rbrace}\abs{\alpha-2u}.
\end{equation}
In order to take advantage of asymptotically flat conditions given in Definition \ref{defn-unif-asym} it we will often consider $\tilde{L}$ sufficiently larger than $\max\lbrace L,R_0\rbrace$. We will then use the above inequality to relate any estimates we obtain back to our original situation. Similarly, we will look at $\tilde{\rho}_1>\max\lbrace\rho_1,R_0\rbrace$.

If we write the area enlarging condition (\ref{defn-area-enl}) in terms of the coordinate functions, then we see that
\begin{equation}
	(\alpha-2u)(\rho_0,z)\ge 0.
\end{equation}
From this, it quickly follows that
\begin{equation}
\int_{\partial\Omega^{\tilde{L}}_{\rho_0\tilde{\rho}_1}\cap\lbrace\rho=\rho_0\rbrace}\abs{\alpha-2u}=\int_{\partial\Omega^{\tilde{L}}_{\rho_0,\tilde{\rho}_1}\cap\lbrace\rho=\rho_0\rbrace}\alpha-2u.
\end{equation}

In order to estimate the above, we once again take advantage of the fundamental theorem of calculus to write
\begin{equation}\label{boundary formula 2u -alpha}
\int_{\partial\Omega^{\tilde{L}}_{\rho_0\tilde{\rho}_1}\cap\lbrace\rho=\rho_0\rbrace}(\alpha-2u) dz=\int_{-\frac{\tilde{L}}{2}}^{\frac{\tilde{L}}{2}}\int_{\rho_0}^{\tilde{\rho}_1}-\frac{\partial(\alpha-2u)}{\partial\rho}d\rho dz+\int_{-\frac{\tilde{L}}{2}}^{\frac{\tilde{L}}{2}}(\alpha-2u)(\tilde{\rho}_1,z)dz.
\end{equation}
We may switch the order of integration for the integral on the right to get
\begin{equation}
\int_{\rho_0}^{\tilde{\rho}_1}\int_{-\frac{\tilde{L}}{2}}^{\frac{\tilde{L}}{2}}-\frac{\partial(\alpha-2u)}{\partial\rho}dzd\rho.
\end{equation}
As before (\ref{normal derivative to a rectangle}), from Stokes' theorem we get
\begin{equation}
\int_{-\frac{\tilde{L}}{2}}^{\frac{\tilde{L}}{2}}-\frac{\partial(\alpha-2u)}{\partial\rho}(\rho,z)dz=\int_{\lbrace\rho\le s,\abs{z}\le\frac{\tilde{L}}{2}\rbrace}\Delta(\alpha-2u)(s,z)-\int_{\lbrace\rho\le s,\abs{z}=\frac{\tilde{L}}{2}\rbrace}\frac{\partial(\alpha-2u)}{\partial\nu}.
\end{equation}
Taking the absolute value of the above and plugging it into (\ref{boundary formula 2u -alpha}) gives us the estimate
\begin{equation}
\int_{-\frac{\tilde{L}}{2}}^{\frac{\tilde{L}}{2}}\abs{\alpha-2u}\le\int_{\rho_0}^{\tilde{\rho}_1}(\int_{\lbrace\rho\le s,\abs{z}=\frac{\tilde{L}}{2}\rbrace}\abs{\Delta(\alpha-2u)}+\int_{\lbrace\rho\le s,\abs{z}=\frac{\tilde{L}}{2}\rbrace}\abs{\frac{\partial(\alpha-2u)}{\partial z}}ds)d\rho+\int_{-\frac{\tilde{L}}{2}}^{\frac{\tilde{L}}{2}}\abs{\alpha-2u}(\tilde{\rho}_1,z)dz.
\end{equation}
We now proceed to estimate the right hand side term by term.

We start with the term
\begin{equation}
\int_{\rho_0}^{\tilde{\rho}_1}\int_{\lbrace\rho\le s,\abs{z}=\frac{\tilde{L}}{2}\rbrace}\abs{\frac{\partial(\alpha-2u)}{\partial z}}dsd\rho.
\end{equation}
Using the asymptotic flatness condition, we estimate
\begin{equation}
\int_{\lbrace\rho\le s,\abs{z}=\frac{\tilde{L}}{2}\rbrace}\abs{\frac{\partial(\alpha-2u)}{\partial z}}ds\le\int_{\lbrace\rho\le s,\abs{z}=\frac{\tilde{L}}{2}\rbrace}\frac{3C}{\abs{(s,z)}^2}ds.
\end{equation}
Once more, a simple integration bounds the above by
\begin{equation}
\frac{6\pi C}{\tilde{L}}.
\end{equation}
Thus, we see that
\begin{equation}
\int_{\rho_0}^{\tilde{\rho}_1}\int_{\lbrace\rho\le s,\abs{z}=\frac{\tilde{L}}{2}\rbrace}\abs{\frac{\partial(\alpha-2u)}{\partial z}}dsd\rho\le\frac{6\pi C\tilde{\rho}_1}{\tilde{L}}.
\end{equation}
We may bound
\begin{equation}
\int_{\rho_0}^{\tilde{\rho}_1}\left(\int_{\lbrace\rho\le s,\abs{z}=\frac{\tilde{L}}{2}\rbrace}\abs{\Delta(\alpha-2u)}\right)d\rho
\end{equation}
by modifying Lemma \ref{L1 Estimate for the Laplacian of 2u-alpha} slightly to get
\begin{equation}
\int_{\lbrace\rho\le s,\abs{z}=\frac{\tilde{L}}{2}\rbrace}\abs{\Delta(\alpha-2u)}\le\frac{4m+4\sqrt{\tilde{L}m}}{\rho}
\end{equation}
and then integrating. We see that
\begin{equation}
\int_{\rho_0}^{\tilde{\rho}_1}\left(\int_{\lbrace\rho\le s,\abs{z}=\frac{\tilde{L}}{2}\rbrace}\abs{\Delta(\alpha-2u)}\right)d\rho\le(4m+4\sqrt{\tilde{L}m})\log(\frac{\tilde{\rho}_1}{\rho_0}).
\end{equation}
Finally, we must bound
\begin{equation}
\int_{-\frac{\tilde{L}}{2}}^{\frac{\tilde{L}}{2}}\abs{\alpha-2u}(\tilde{\rho}_1,z)dz.
\end{equation}
Oddly enough, this turns out to be the most delicate estimate, and the point where we need our extra assumption on the asymptotic falloff of the function $\alpha$. From Lemma \ref{Uniform Control of u}, we know that the $C^{0,\beta}$ norm of $u$ is controlled by $m$. Recalling (\ref{uniform estimate for u}), we see that there is a constant $\tilde{\epsilon}(\tilde{\rho}_1,m)$ such that
\begin{equation}
\int_{-\frac{\tilde{L}}{2}}^{\frac{\tilde{L}}{2}}\abs{u(\tilde{\rho}_1,z)}dz\le \tilde{L}\tilde{\epsilon}(m,\tilde{\rho}_1).
\end{equation}
Again, looking at Lemma \ref{Uniform Control of u}, we see that for fixed $\tilde{\rho}_1$
\begin{equation}
	\lim_{m\rightarrow 0}\tilde{\epsilon}(\tilde{\rho}_1,m)=0.
\end{equation}
From the extra assumption on the asymptotic falloff of $\alpha$, we see that
\begin{equation}
\int_{-\frac{\tilde{L}}{2}}^{\frac{\tilde{L}}{2}}\abs{\alpha(\tilde{\rho}_1,z)}dz\le\int_{-\frac{\tilde{L}}{2}}^{\frac{\tilde{L}}{2}}\frac{C}{\abs{(\tilde{\rho}_1,z)}^{1+\tau}}dz\le C(\tau)(\tilde{\rho}_1)^{-\tau},
\end{equation}
where $C(\tau)$ is a constant depending only on $\tau$.
We may put all of this together to see that
\begin{equation}\label{area enlarging interior boundary estimate alpha-2u}
\int_{-\frac{\tilde{L}}{2}}^{\frac{\tilde{L}}{2}}\abs{\alpha-2u}dz\le(4m+4\sqrt{\tilde{L}m})\log(\frac{\tilde{\rho}_1}{\rho_0})+\frac{6\pi C\tilde{\rho}_1}{\tilde{L}}+\tilde{L}\tilde{\epsilon}(\tilde{\rho}_1,m)+C(\tau)(\tilde{\rho}_1)^{-\tau}.
\end{equation}
By choosing $\tilde{\rho}_1$ and $\tilde{L}$ to be as large as necessary and choosing $m$ to be as small as necessary, we see that the above quantity can be made as small as we desire. $\square$

The following corollary to Lemma \ref{2u-alpha boundary estimate} is analogous to Corollary \ref{Boundary Estimate for alpha-2u}.
\begin{cor}\label{boundary control area enlarging case}
	Let $\mathcal{M}$ be a family of axisymmetric metrics with nonnegative scalar curvature which is strongly uniformly asymptotically flat outside of radius $R_0$. Suppose also that $M$ is area enlarging at $\rho_0$. Let $\Omega^{L}_{\rho_0\rho_1}$ denote the region $\lbrace(\rho,z)|\rho_0\le\rho\le \rho_1,\abs{z}\le\frac{L}{2}\rbrace$, and $(\Omega^{L}_{\rho_0\rho_1})_{\sigma}$ denote $\lbrace x\in\Omega^{L}_{\rho_0\rho_1}|d(x,\partial\Omega^{\rho_1}_{\rho_0})>\sigma\rbrace$. Then for $m>0$, $\sigma>0$, $L>R_0$, and $\rho_1>R_)$ there is a constant $C(\tau,m,\sigma,L,\rho_1,\rho_0)$ such that if $g\in\mathcal{M}$ and the ADM mass of $g$ is less than $m$,
	then
	\begin{equation}\label{sup norm of boundary}
		\sup_{x\in(\Omega^{L}_{\rho_0\rho_1})_\sigma}\exp\left(\int_{\partial\Omega^{L}_{\rho_0\rho_1}}|H_D(x,y)\frac{\partial(\alpha-2u)}{\partial\nu}(y)|+|(\alpha-2u)(y)\frac{\partial H_D}{\partial\nu}(x,y)|dy\right)\le \exp\left[C(\tau,m,\sigma,L,\rho_1,\rho_0)\right]
	\end{equation}
where $\tau$ is the constant appearing in (\ref{Str-Unf-Asym-Flat}) and $C(,\tau,m,\sigma,L,\rho_1,\rho_0)$ is a constant depending on $\tau$, $m$, $\sigma$, $L$, $\rho_1$, and $\rho_0$. 
\end{cor}
Proof: Much of the proof remains the same as it was in the radially monotone case. The only difference is that we need to estimate
\begin{equation}
\int_{\partial\Omega^{L}_{\rho_0\rho_1}\cap\lbrace\rho=\rho_0\rbrace}\abs{\alpha-2u},
\end{equation}
instead of 
\begin{equation}
\int_{\partial\Omega^{L}_{\rho_0\rho_1}\cap\lbrace\rho=\rho_0\rbrace}\abs{\frac{\partial(\alpha-2u)}{\partial\nu}}.
\end{equation}
This we did in Lemma \ref{2u-alpha boundary estimate}. $\qed$

We now estimate the $W^{1,p}$ norm of $\alpha-2u$. Using the function $H_D$ instead of $H_N$ complicates our estimate of $\abs{\abs{\nabla(\alpha-2u)}}_{L^p\left(\Omega^{L}_{\rho_0\rho_1}\right)}$.  We resort to shrinking our region a bit.
\begin{lemma}\label{sobolev estimates alpha-2u area enlarging case}
	Let $\mathcal{M}$ be a family of axisymmetric metrics with nonnegative scalar curvature which is strongly uniformly asymptotically flat outside of radius $R_0$. Suppose also that $\mathcal{M}$ is area enlarging at $\rho_0$. For any $\rho_1>\rho_0$, $L$, $1\le p<2$, $\sigma>0$, and $\epsilon>0$ there is a $\delta>0$ such that if $g\in\mathcal{M}$ and
	\begin{equation}
		m(g)<\delta,
	\end{equation}
	then
	\begin{equation}
		\abs{\abs{\alpha-2u}}_{W^{1,p}\left(\left(\Omega^{L}_{\rho_0\rho_1}\right)_{\sigma}\right)}<\epsilon.
	\end{equation}
	Here
	\begin{equation}
		\left(\Omega^{L}_{\rho_0\rho_1}\right)_{\sigma}:=\lbrace x\in\Omega^{L}_{\rho_0\rho_1}:d(x,\partial\Omega^{L}_{\rho_0\rho_1})\ge\sigma\rbrace.
	\end{equation}
\end{lemma}
\textbf{Proof}: We may estimate the $L^p$ norm of $\alpha-2u$ much as we did in Lemma (\ref{Sobolev Convergence of 2u-alpha}). We once again consider $\tilde{L}>L$ and $\tilde{\rho}_1>\rho_0$. As before,
\begin{equation}
	\int_{\left(\Omega^{L}_{\rho_0\rho_1}\right)_{\sigma}}\abs{\alpha-2u}^p\le C(p)\int_{\left(\Omega^{L}_{\rho_0\rho_1}\right)_{\sigma}}\left(\int_{\Omega^{\tilde{L}}_{\rho_0\tilde{\rho}_1}}\abs{(\alpha-2u)\frac{\partial H_D}{\partial\nu}}+\abs{H_D\frac{\partial(\alpha-2u)}{\partial\nu}}\right)^p+\left(\int_{\Omega^{\tilde{L}}_{\rho_0\tilde{\rho}_1}}\abs{H_D\Delta(\alpha-2u)}\right)^pdx.
\end{equation}
On $\partial\Omega^{\tilde{L}}_{\rho_0\tilde{\rho}_1}-\lbrace\rho=\rho_0\rbrace$ we have the following bound on the boundary terms
\begin{equation}
	\frac{24C\tilde{\rho}_1}{\pi\tilde{L}\abs{\tilde{L}-L}}+\frac{3C\tilde{L}}{\pi\tilde{\rho}_1\abs{\tilde{\rho}_1-\rho_1}}+\frac{24C\tilde{\rho}_1\log\left(2\sqrt{\tilde{L}^2+\tilde{\rho}^2_1}\right)}{\pi\tilde{L}^2}+\frac{3C\tilde{L}\log\left(\sqrt{\tilde{L}^2+\tilde{\rho}^2_1}\right)}{\pi\tilde{\rho}^2_1}.
\end{equation}
Using the proof of Lemma \ref{2u-alpha boundary estimate} for terms on $\partial\Omega^{\tilde{L}}_{\rho_0\tilde{\rho}_1}\cap\lbrace\rho=\rho_0\rbrace$, we have the estimate
\begin{equation}
	\frac{1}{\pi\sigma}\left((4m+4\sqrt{\tilde{L}m})\log(\frac{\tilde{\rho}_1}{\rho_0})+\frac{6\pi C\tilde{\rho}_1}{\tilde{L}}+\tilde{L}\tilde{\epsilon}(\tilde{\rho}_1,m)+C(\tau)(\tilde{\rho}_1)^{-\tau}\right).
\end{equation}
If we let $\tilde{\rho}_1=\tilde{L}^{\frac{2}{3}}$, then we may see that we may pick $\tilde{L}$ large enough and $m$ small enough to ensure
\begin{equation}
	\abs{\abs{\alpha-2u}}_{L^p\left(\left(\Omega^{L}_{\rho_0\rho_q}\right)_\sigma\right)}<\frac{\epsilon}{2}.
\end{equation}	

If we differentiate Green's representation formula with $H_D$ we get
\begin{equation}
\nabla(\alpha-2u)(x)=\int_{\partial\Omega^{\tilde{L}}_{\rho_0\tilde{\rho}_1}}(\alpha-2u)\nabla_x\left(\frac{\partial H_D}{\partial\nu}\right)-\nabla_x\left(H_D(x,y)\right)\frac{\partial(\alpha-2u)}{\partial\nu}dy+\int_{\Omega^{\tilde{L}}_{\rho_0\tilde{\rho}_1}}\nabla_x\left(H_D(x,y)\right)\Delta(\alpha-2u)dy.
\end{equation}
On $\partial\Omega^{\tilde{L}}_{\rho_0\tilde{\rho}_1}\cap\lbrace\rho=\rho_0\rbrace$ the above expression is particularly difficult to work with. The issue is that we cannot integrate
\begin{equation}
\abs{\nabla_x\left(\frac{\partial H_D}{\partial\nu}\right)}\sim\frac{1}{\abs{x-y}^2}
\end{equation}
for $x$ near the boundary, and so we cannot complete the estimate of $\abs{\abs{\alpha-2u}}_{W^{1,p}}$ in the same way we proved  \ref{Sobolev Convergence of 2u-alpha}.

As we have done before, we take the absolute value of both sides and raise the result to the power $p$ and then integrate to see that
\begin{equation}
	\int_{\left(\Omega^{\tilde{L}}_{\rho_0\tilde{\rho}_1}\right)_{\sigma}}\abs{\nabla(\alpha-2u)}^p
\end{equation}
is bounded above by
\begin{equation}\label{area enlarging gradient estimate of alpha-2u}
	 C(p)\int_{\left(\Omega^{L}_{\rho_0\rho_1}\right)_{\sigma}}\left(\int_{\partial\Omega^{\tilde{L}}_{\rho_0\tilde{\rho}_1}}\abs{\frac{\partial(\alpha-2u)}{\partial\nu}\nabla_xH_D}+\abs{(\alpha-2u)\nabla_x\frac{\partial H_D}{\partial\nu}}dy\right)^p+\left(\int_{\Omega^{\tilde{L}}_{\rho_0\tilde{\rho}_1}}\abs{\Delta(\alpha-2u)\nabla_xH_D}dy\right)^pdx.
\end{equation}
We  once again split the first term into the following two pieces:
\begin{equation}
	\partial\Omega^{\tilde{L}}_{\rho_0\tilde{\rho}_1}-\lbrace\rho=\rho_0\rbrace
\end{equation}
and
\begin{equation}
	\partial\Omega^{\tilde{L}}_{\rho_0\tilde{\rho}_1}\cap\lbrace\rho=\rho_0\rbrace.
\end{equation}
Both pieces are relatively easy to estimate. For the first piece the estimates are similar to the above.

As was noted earlier, the gradient of $\nabla_x\frac{\partial H_D}{\partial\nu}$ isn't integrable over $\Omega^{L}_{\rho_0\rho_1}$ for $y$ in $\partial\Omega^{\tilde{L}}_{\rho_0\tilde{\rho}_1}\cap\lbrace\rho=\rho_0\rbrace$. However, $\nabla_x\frac{\partial H_D}{\partial\nu}$ is much better behaved away from $\partial\Omega^{L}_{\rho_0\rho_1}$. We now attempt to estimate
\begin{equation}
	\int_{\left(\Omega^{L}_{\rho_0\rho_1}\right)_{\sigma}}\left(\int_{\partial\Omega^{\tilde{L}}_{\rho_0\tilde{\rho}_1}}\abs{(\alpha-2u)\nabla_x\frac{\partial H_D}{\partial\nu}}dy\right)^pdx.
\end{equation}
As we did before, we split $\partial\Omega^{\tilde{L}}_{\rho_0\tilde{\rho}_1}\cap\lbrace\rho=\rho_0\rbrace$ into
\begin{equation}\label{Inner Piece One Area-Enlg Case}
	\partial\Omega^{\tilde{L}}_{\rho_0\tilde{\rho}_1}\cap\lbrace\rho=\rho_0,\abs{z}\le L\rbrace
\end{equation}
and
\begin{equation}\label{Inner Piece Two Area-Enlg Case}
	\partial\Omega^{\tilde{L}}_{\rho_0\tilde{\rho}_1}\cap\lbrace\rho=\rho_0,\abs{z}>L\rbrace.
\end{equation}
We start with the piece (\ref{Inner Piece One Area-Enlg Case}).
We may use Minkowski's integral inequality \cite{Folland} to see that
\begin{equation}
	\left(\int_{\left(\Omega^{L}_{\rho_0\rho_1}\right)_{\sigma}}\left(\int_{\partial\Omega^{\tilde{L}}_{\rho_0\tilde{\rho}_1}\cap\lbrace\rho=\rho_0,\abs{z}\le L\rbrace}\abs{(\alpha-2u)\nabla\frac{\partial H_D}{\partial\nu}}dy\right)^pdx\right)^{\frac{1}{p}}
\end{equation}
is bounded above by
\begin{equation}
	\int_{\partial\Omega^{\tilde{L}}_{\rho_0\tilde{\rho}_1}\cap\lbrace\rho=\rho_0,\abs{z}\le L\rbrace}\abs{\alpha-2u}\left(\int_{\left(\Omega^{L}_{\rho_0\rho_1}\right)_{\sigma}}\abs{\nabla\frac{\partial H_D}{\partial\nu}}^pdx\right)^{\frac{1}{p}}dy.
\end{equation}
We now estimate 
\begin{equation}
	\int_{\left(\Omega^{L}_{\rho_0\rho_1}\right)_{\sigma}}\abs{\nabla\frac{\partial H_D}{\partial\nu}}^pdx
\end{equation}
for $y$ in $\partial\Omega^{\tilde{L}}_{\rho_0\tilde{\rho}_1}\cap\lbrace\rho=\rho_0,\abs{z}\le L\rbrace$. Both $\partial\Omega^{\tilde{L}}_{\rho_0\tilde{\rho}_1}\cap\lbrace\rho=\rho_0,\abs{z}\le L\rbrace$ and $\left(\Omega^{L}_{\rho_0\rho_1}\right)_{\sigma}$ are contained in $\Omega^{2L}_{\rho_0\rho_1}$. Thus, if we let $r_0$ be the diameter of $\Omega^{2L}_{\rho_0\rho_1}$, then for all $y\in\partial\Omega^{\tilde{L}}_{\rho_0\tilde{\rho}_1}\cap\lbrace\rho=\rho_0,\abs{z}\le L\rbrace$ we have
\begin{equation}
	\int_{\left(\Omega^{L}_{\rho_0\rho_1}\right)_{\sigma}}\abs{\nabla_x\frac{\partial H_D}{\partial\nu}}^p\le\int_{B(y,r_0)\backslash B(y,\sigma)}\frac{3^p}{\pi^p\abs{x-y}^{2p}}dx=3^p\pi^{1-p}2 \int_{\sigma}^{r_0}Cr^{-2p+1}dr=C(p,L,\rho_1,\sigma).
\end{equation}
Thus, we may see that
\begin{equation}\label{interior gradient boundary estimate piece one}
	\left(\int_{\left(\Omega^{L}_{\rho_0\rho_1}\right)_{\sigma}}\left(\int_{\partial\Omega^{\tilde{L}}_{\rho_0\tilde{\rho}_1}\cap\lbrace\rho=\rho_0,\abs{z}\le L\rbrace}\abs{(\alpha-2u)\nabla\frac{\partial H_D}{\partial\nu}}dy\right)^pdx\right)^{\frac{1}{p}}\le C(p,L,\rho_1,\sigma)^{\frac{1}{p}}\int_{\partial\Omega^{\tilde{L}}_{\rho_0\tilde{\rho}_1}\cap\lbrace\rho=\rho_0,\abs{z}\le L\rbrace}\abs{\alpha-2u}.
\end{equation}

Over (\ref{Inner Piece Two Area-Enlg Case}) we have
\begin{equation}
	\abs{\nabla\frac{H_D}{\partial\nu}}\le\frac{12}{\pi L^2}.
\end{equation}
Thus, we have
\begin{equation}\label{interior gradient boundary estimate piece two}
	\left(\int_{\left(\Omega^{L}_{\rho_0\rho_1}\right)_{\sigma}}\left(\int_{\partial\Omega^{\tilde{L}}_{\rho_0\tilde{\rho}_1}\cap\lbrace\rho=\rho_0,\abs{z}>L\rbrace}\abs{(\alpha-2u)\nabla\frac{\partial H_D}{\partial\nu}}dy\right)^pdx\right)^{\frac{1}{p}}\le\left(\frac{12\rho_1}{L}\right)^{\frac{1}{p}}\int_{\partial\Omega^{\tilde{L}}_{\rho_0\tilde{\rho}_1}\cap\lbrace\rho=\rho_0,\abs{z}>L\rbrace}\abs{(\alpha-2u)}dy.
\end{equation}
For the last term in (\ref{area enlarging gradient estimate of alpha-2u}) we may use the Riesz potential estimate as we have done before.
Putting everything together gives us the result. $\qed$

In fact, the steps required in the above proof give us a corollary analogous to \ref{2u-alpha Blowup  Estimate}.
\begin{cor}\label{2u-alpha Blowup Estimate Area-Enlg Case}
	Let $\mathcal{M}$ be a family of axisymmetric metrics with nonnegative scalar curvature which is strongly uniformly asymptotically flat outside of radius $R_0$. Suppose $\mathcal{M}$ is area enlarging as well. For any $L$, $\rho_1$, $1\le p<2$, and $\epsilon>0$ there exist a $\delta>0$ such that if $g\in\mathcal{M}$ and
	\begin{equation}
		m(g)<\delta,
	\end{equation}
	then
	\begin{equation}\label{area enlarging controlling blowup}
		\int_{\Omega^{L}_{\rho_0,\rho_1}}\abs{\alpha-2u}^p<\frac{\epsilon\abs{\log\rho_0}^p}{\rho_0^p}
	\end{equation}
	and
	\begin{equation}\label{area enlarging controlling gradient blowup}
		\int_{\Omega^{L}_{\rho_0,\rho_1}}\abs{\nabla(\alpha-2u)}^p\le\frac{\epsilon\abs{\log\rho_0}^p}{\rho_0^p}.
	\end{equation}
\end{cor}
\textbf{Proof}: The proofs of (\ref{area enlarging controlling blowup}) and (\ref{area enlarging controlling gradient blowup}) are similar. We only prove (\ref{area enlarging controlling gradient blowup}). Observe that
\begin{equation}
	\Omega^{L}_{2\rho_0\rho_1}\subset\left(\Omega^{L+\sigma}_{\rho_0(\rho_1+\sigma)}\right)_{\sigma}.
\end{equation}
In particular, we see from the estimates in the above theorem that
\begin{equation}
	\int_{\Omega^{L}_{2\rho_0\rho_1}}\abs{\nabla(\alpha-2u)}^p\le\int_{\left(\Omega^{L+\sigma}_{\rho_0(\rho_1+\sigma)}\right)_{\sigma}}\abs{\nabla(\alpha-2u)}^p
\end{equation}
is bounded above by
\begin{equation}\label{Blowup estimate intermediate}
	C(p,L,\rho_1,\sigma)\left[(4m+4\sqrt{\tilde{L}m})\log(\frac{\tilde{\rho}_1}{\rho_0})+D(m,\tilde{L},\tilde{\rho}_1,\tau)\right]^p+E(p,\tilde{L},\tilde{\rho}_1)\left(\frac{4m+4\sqrt{\tilde{L}m}}{\rho_0}\right)^p+
	F(p,\tilde{L},\tilde{\rho}_1),
\end{equation}
where $C(p,L,\rho_1,\sigma)$ is a combination of the constants found in (\ref{interior gradient boundary estimate piece one}) and (\ref{interior gradient boundary estimate piece two}), $D(m,\tilde{L},\tilde{\rho}_1,\tau)$ is the remainder of (\ref{area enlarging interior boundary estimate alpha-2u}),  $E(p,\tilde{L},\tilde{\rho}_1)$ comes from the Riesz potential estimate, and $F(p,\tilde{L},\tilde{\rho}_1)$ is the bound on the remaining boundary terms estimated in (\ref{area enlarging gradient estimate of alpha-2u}).
A simple calculation shows that for $1<p<2$
\begin{equation}
	C(L,\rho_1,\sigma)\le C(p)\sigma^{-p},
\end{equation}
since $2-2p>-p$.
For $p=1$, we have
\begin{equation}
	C(L,\rho_1,\sigma)\le C(L,\rho_1)\log(\sigma).
\end{equation}
If we plug the above into (\ref{Blowup estimate intermediate}) with $\sigma=\rho_0$, then we may see that choosing $\tilde{L}$ and $\tilde{\rho}_1$ large enough, and choosing mass to be small enough gives the result. $\qed$

We may now prove a theorem analogous to Theorem \ref{Sobolev Convergence of 2u-alpha}.
\begin{lemma}\label{Area Assumption Sobolev Estimate for 2u-alpha}
	Let $\mathcal{M}$ be an uniformly asymptotically flat family of metrics with nonnegative scalar curvature. Suppose that $M$ be area enlarging. Let $\Omega^{L}_{\rho_0\rho_1}$ denote the rectangle given by $\lbrace(\rho,z)|\rho_0\le\rho\le \rho_1,\abs{z}\le\frac{L}{2}\rbrace$ and let $(\Omega^{L}_{\rho_0\rho_1})_{\sigma}$ denote $\lbrace x\in\Omega^{L}_{\rho_0\rho_1}|d(x,\partial\Omega^{L}_{\rho_0\rho_1})>\sigma\rbrace$. For any $1\le p<2$, $\sigma>0$, $\rho_0>0$, and $\epsilon>0$ there exists a $\delta>0$ such that if $g$ is in our collection of uniformly asymptotically flat metrics, the ADM mass of $g$ is less than $\delta$, and, in the axisymmetric coordinate representation of $g$
	then
	\begin{equation}
	\abs{\abs{e^{\abs{\alpha-2u}}-1}}_{W^{1,p}((\Omega^{L}_{\rho_0\rho_1})_\sigma)}<\epsilon.
	\end{equation}
\end{lemma}

\textbf{Proof}: The proof follows the same line as in the radially monotone case,
except we use Lemma \ref{sobolev estimates alpha-2u area enlarging case} instead of Proposition \ref{Sobolev Convergence of 2u-alpha}. It can be shown that Corollary \ref{Msr-like corollary} can be adapted to the function $H_D$. Thus, we also use Corollary \ref{boundary control area enlarging case} instead of Lemma \ref{Boundary Estimate for alpha-2u}. $\qed$

Now that we have analogues of all the estimates we made in the radially monotone case, the proofs of Theorem \ref{W1p estimates}, Theorem \ref{Volume Estimates}, Theorem \ref{Area Estimates}, Theorem \ref{d estimates} and Theorem \ref{ufm con. in the ext.} follow almost exactly as they did in the radially monotone case. The only theorem whose modification to the area-enlarging case requires a little care is Theorem \ref{ufm con. in the ext.}. Since Corollary \ref{2u-alpha Blowup Estimate Area-Enlg Case} has a slightly different hypothesis than Corollary \ref{2u-alpha Blowup  Estimate}, we must show that the conclusion of Lemma \ref{Countering Blowup Estimate} holds with a slightly weaker hypothesis.
\begin{lemma}
	Let $f$ be a measurable function on $\Omega^{L}_{0\rho_1}$. Suppose for each $t$ we have the estimate
	\begin{equation}
	\int_{\Omega^{\rho_1}_{t}}\abs{f}\le\frac{\epsilon\abs{\log(t)}^{\tilde{q}}}{t^{q}}
	\end{equation}
	for some $\epsilon>0$, $q>0$, and $\tilde{q}$. Suppose $\sigma>q$. Then, there exists a constant, denoted $C(\sigma,q,\tilde{q})$, depending only on $\sigma$, $q$, and $\tilde{q}$ such that
	\begin{equation}
	\int_{\Omega^{\rho_1}_{0}}\rho^{\sigma}\abs{f}\le C(\sigma,q,\tilde{q})\epsilon.
	\end{equation}
\end{lemma}
\textbf{Proof}: As before, let $t_n=2^{-n}\rho_1$ and 
let $\Omega_{t_n,t_{n-1}}$ be the following rectangle.
\begin{equation}
\Omega_{t_n,t_{n-1}}=\lbrace t_n\le\rho\le t_{n-1},\abs{z}\le\frac{L}{2}\rbrace
\end{equation}
From the Monotone Convergence Theorem we see that
\begin{equation}
\int_{\Omega_{0,t_0}}\rho^{\sigma}\abs{f}=\sum_{1}^{\infty}\int_{\Omega_{t_n,t_{n-1}}}\rho^{\sigma}\abs{f}.
\end{equation}
We now make the estimate
\begin{equation}
\int_{\Omega_{t_n,t_{n-1}}}\rho^{\sigma}\abs{f}\le t^{\sigma}_{n-1}\frac{\epsilon\abs{\log(t_n)}^{\tilde{q}}}{t^{q}_{n}}=2^{q}\rho^{\sigma-q}_1(2^{\sigma-q})^{-n}\abs{\log(2^{-n}\rho_1)}^{\tilde{q}}\epsilon.
\end{equation}
This gives a convergent series so long as $\sigma>q$, where we have used that $\sigma-q=\lambda>0$ and
\begin{equation}
	\lim_{n\rightarrow\infty}\rho_12^{-n}\abs{\log(\rho_12^{-n})}^{\frac{2\tilde{q}}{\lambda}}=0.
\end{equation}
In total, we have the estimate
\begin{equation}
\int_{\Omega_{0,t_0}}\rho^{\sigma}\abs{f}\le C(\sigma,q,\tilde{q})\epsilon.\qed
\end{equation}

Now we can show that Lemma \ref{2u-alpha Sobolev Estimate over B_R} holds in the area-enlarging case and so Theorem \ref{ufm con. in the ext.} also holds in the area-enlarging case.
\appendix
\section[short title]{Examples}
\label{secA} \setcounter{equation}{0}
\setcounter{section}{1}
\subsection{Kerr-Newman}
In this section, we show that the Kerr-Newman family of metrics satisfy both the radial monotone condition and the area enlarging condition. This is done by a direct calculation. We take the familiar Brill-Lindquist coordinates and transform them into cylindrical coordinates. Unfortunately, the simple expression of the Kerr-Newman metric in Brill-Lindquist coordinates becomes rather complicated when it is written in cylindrical coordinates. The procedure itself is uncomplicated, since there is an explicit map between these two coordinates. The change of coordinates depends on the charge, angular momentum, and mass of the Kerr-Newman metric. Once the map has been constructed, we use the expression for the metric in Brill-Lindquist to write down the expression for the metric in cylindrical coordinates.

We now describe in detail the coordinate change from Brill-Lindquist coordinates to cylindrical coordinates and write down the exact formula for the metric functions $u$ and $\alpha$. It is convenient to introduce a third coordinate system between Brill-Lindquist and cylindrical. We shall use the Prolate-Spheroidal coordinates. We will first consider the map from Prolate Spheroidal coordinates to Brill-Lindquist coordinates, and then pull back the metric. Let $a$ denote the angular momentum parameter, let $e$ denote the charge parameter, and let $m$ denote the mass parameter, then, in Brill-Lindquist coordinates, the Kerr metric takes the form
\begin{equation}
g=\frac{\sigma}{\gamma}dr^2+\sigma d\theta^2+\frac{\sin^2(\theta)}{\sigma}[(r^2+a^2)^2-a^2\sin^2(\theta)\gamma(r)]d\phi^2
\end{equation}
for
\begin{equation}
\gamma(r)=r^2-2mr+a^2+e^2
\end{equation}
and
\begin{equation}
\sigma(r,\theta)=r^2+a^2\cos^2(\theta).
\end{equation}

The map from prolate spheroidal coordinates $(x,y,\phi)$ to Brill-Lindquist coordinates $(r,\theta,\phi)$ is given by
\begin{equation}
r=x\sqrt{m^2-(a^2+e^2)}+m
\end{equation}
\begin{equation}
\theta=\cos^{-1}(y)
\end{equation}
For convenience, we will write
\begin{equation}
k=\sqrt{m^2-(a^2+e^2)}.
\end{equation}
The map from cylindrical coordinates to prolate spheroidal is, unfortunately, less simple.
\begin{equation}
x=\frac{\sqrt{\rho^2+(z+k)^2}+\sqrt{\rho^2+(z-k)^2}}{2k}
\end{equation}
\begin{equation}
y=\frac{\sqrt{\rho^2+(z+k)^2}-\sqrt{\rho^2+(z-k)^2}}{2k}
\end{equation}

We now pull back the Kerr-Newman metric twice to obtain the formulas for the functions $u$ and $\alpha$ in cylindrical coordinates. The end results of this process are the following formulas
\begin{equation}
u(\rho,z)=-\frac{1}{2}\log[\frac{(1-y^2)([(kx+m)^2+a^2]^2-a^2k^2[1-y^2][x^2-1])}{\rho^2([kx+m]^2+a^2y^2)}]
\end{equation}
\begin{equation}
\alpha(\rho,z)=\frac{1}{2}\log[\frac{(kx+m)^2+a^2y^2}{k^2(x^2-y^2)}]+u(\rho,z)
\end{equation}
When written entirely in terms of $(\rho,z)$, these two equations are very cumbersome. Luckily, for the purpose of verifying the radial monotonicity condition and the area enlarging condition, writing everything in terms of $(\rho,z)$ turns out to be unnecessary.

A straight forward calculation shows that
\begin{equation}\label{rho derivative in prolate spheroidal}
	\frac{\partial}{\partial\rho}=\frac{\rho}{\left(\rho^2+(z+k)^2\right)^{\frac{1}{2}}\left(\rho^2+(z-k)^2\right)^\frac{1}{2}}\left(x\frac{\partial}{\partial x}-y\frac{\partial}{\partial y}\right).
\end{equation}
Thus, we see that
\begin{equation}
	\frac{\partial(\alpha-2u)}{\partial\rho}=f(\rho,z)\left(x\frac{\partial}{\partial x}-y\frac{\partial}{\partial y}\right)\log\left(\frac{\left[(kx+m)^2+a^2\right]^2-a^2k^2[1-y^2][x^2-1]}{k^2(x^2-1)(x^2-y^2)}\right),
\end{equation}
where $f(\rho,z)$ is the nonnegative function appearing in front of the derivatives in (\ref{rho derivative in prolate spheroidal}). Since $f(\rho,z)$ is nonnegative, we may restrict our analysis to the second term on the right. Taking the derivatives and collecting terms leaves us with
\begin{equation}
	\begin{split}
	& \frac{4kx(kx+m)\left[(kx+m)^2+a^2\right]-2a^2k^2x^2(1-y^2)}{\left[(kx+m)^2+a^2\right]^2-a^2k^2(1-y^2)(x^2-1)}-\frac{2x^2\left((x^2-1)+(x^2-y^2)\right)}{(x^2-1)(x^2-y^2)}+\\ &-\left[\frac{2a^2k^2(x^2-1)y^2}{\left[(kx+m)^2+a^2\right]^2-a^2k^2(1-y^2)(x^2-1)}+\frac{2y^2}{x^2-y^2}\right]
	\end{split}.
\end{equation}
The third term in brackets is nonnegative, so we must analyze the interplay of the first two terms.

We expand
\begin{equation}
	\frac{2x^2\left((x^2-1)+(x^2-y^2)\right)}{(x^2-1)(x^2-y^2)}
\end{equation}
to
\begin{equation}
	\frac{2x^2}{x^2-1}+\frac{2x^2}{x^2-y^2}.
\end{equation}
From the range of values that $x$ and $y$ can take, we may deduce that the denominators of both fractions are smaller than $x^2$. Thus, we have
\begin{equation}
	\frac{2x^2}{x^2-1}+\frac{2x^2}{x^2-y^2}>4.
\end{equation}
We now observe that
\begin{equation}
	\left[(kx+m)^2+a^2\right]^2-a^2k^2(1-y^2)(x^2-1)\ge(kx+m)^4+a^2(kx+m)^2.
\end{equation}
As a consequence, we have that
\begin{equation}
	\frac{4kx(kx+m)\left[(kx+m)^2+a^2\right]-2a^2k^2x^2(1-y^2)}{\left[(kx+m)^2+a^2\right]^2-a^2k^2(1-y^2)(x^2-1)}\le 4.
\end{equation}
Putting everything together shows that
\begin{equation}
	\frac{\partial(\alpha-2u)}{\partial\rho}<0. \qed
\end{equation}

It is interesting to explore some of the geometric meaning behind the condition of radial monotonicity. In coordinates, radial monotonicity implies that
\begin{equation}
	\frac{\partial(\alpha-2u)}{\partial\rho}\le 0.
\end{equation}
Recall from the proof of Proposition \ref{Area Enlarging} that the coordinate function $\alpha-2u$ controls the area of axisymmetric surfaces. Thus, it is reasonable to suppose that the radial monotonicity condition is an assumption on the mean curvature of the level sets of the function $\rho$, which is the solution to (\ref{rho equation}). It turns out that this is the case, although in a slightly round about way.
\begin{prop}\label{Sub-IMCF Radial Monotonicity}
	Suppose that $g$ is an asymptotically flat axisymmetric metric and $\rho$ is the solution to (\ref{rho equation}) for $g$. The metric $g$ is radially monotone if and only if the level sets of $\rho$ form a family of surfaces evolving by a sub-inverse-mean-curvature flow. 
\end{prop}
\textbf{Proof}: Let $\eta$ denote the killing field generating the axisymmetry of $(M,g)$. We start by observing that we may lift any function $\omega$ on $M/S^1$ to a function on $M$, which we also denote $\omega$. When considered as a function on $M$ we have
\begin{equation}
	g(\nabla\omega,\eta)=0,
\end{equation}
since we lifted $\omega$ by transporting it along the flow lines of $\eta$. Let $q$ denote the orbit metric of $M/S^1$.
Recall that
\begin{equation}
	q(X,Y)=g\left(\bar{X},\bar{Y}\right)-\frac{g\left(\bar{X},\eta\right)g\left(\bar{Y},\eta\right)}{|\eta|^2_g},
\end{equation}
where $X$ and $Y$ are the images of $\bar{X}$ and $\bar{Y}$ under the projection map, respectively. From the above, we may conclude that for any two functions $\omega$ and $h$ on $M/S^1$ we have
\begin{equation}
	q(\nabla\omega,\nabla h)=g(\nabla\omega,\nabla h).
\end{equation}
We have abused notation slightly in using $\nabla$ to denote both the gradient in $(M/S^1,q)$ and in $(M,g)$.

It is a standard computation to see that the mean curvature of the level sets of $\rho$ is given by
\begin{equation}\label{mean curvature equation}
	H=div_g\left(\frac{\nabla\rho}{\abs{\nabla\rho}_g}\right).
\end{equation}
We expand out the right hand side to get
\begin{equation}
	div_g\left(\frac{\nabla\rho}{\abs{\nabla\rho}_g}\right)=\frac{1}{\abs{\nabla\rho}_g}\left(\Delta_g\rho-\frac{g\left(\nabla\rho,\nabla\abs{\nabla\rho}\right)}{\abs{\nabla\rho}}\right)
\end{equation}
We now use the equation for $\rho$ (\ref{rho equation}) to rewrite the above as
\begin{equation}\label{divergence of normal to level sets}
	\frac{1}{\abs{\nabla\rho}}\left(\frac{g(\nabla\rho,\nabla\abs{\eta})}{\abs{\eta}}-\frac{g\left(\nabla\rho,\nabla\abs{\nabla\rho}\right)}{\abs{\nabla\rho}}\right)=\frac{1}{\abs{\nabla\rho}}g\left(\nabla\rho,\nabla\log\frac{\abs{\eta}}{\abs{\nabla\rho}}\right).
\end{equation}
From axisymmetry, $\abs{\nabla\rho}$ and $\abs{\eta}$ are functions on $M/S^1$. In particular
\begin{equation}
	g\left(\nabla\rho,\nabla\log\frac{\abs{\eta}}{\abs{\nabla\rho}}\right)=q\left(\nabla\rho,\nabla\log\frac{\abs{\eta}}{\abs{\nabla\rho}}\right).
\end{equation}
Recalling the radial monotonicity condition (\ref{defn-rad-mono}) and noting that $\log$ is a monotone increasing function, we see that
\begin{equation}\label{modified radial monotonicity}
	q\left(\nabla\rho,\nabla\log\left(\frac{\abs{\eta}}{\rho\abs{\nabla\rho}}\right)\right)\le 0,
\end{equation}
since in the orbit space $M/S^1$ we have
\begin{equation}
	\frac{\partial}{\partial\rho}=\abs{\frac{\partial}{\partial\rho}}_q^2\nabla\rho.
\end{equation}
We may plug (\ref{divergence of normal to level sets}) and (\ref{mean curvature equation}) into (\ref{modified radial monotonicity}) to see that
\begin{equation}
	0\ge q\left(\nabla\rho,\nabla\log\left(\frac{\abs{\eta}}{\abs{\nabla\rho}}\right)\right)-q\left(\nabla\rho,\nabla\log\rho\right)=\abs{\nabla\rho}H-\abs{\nabla\rho}\abs{\nabla\log\rho}.
\end{equation}
Dividing both sides by $\abs{\nabla\rho}$ and rearranging terms gives
\begin{equation}
	\abs{\nabla\log\rho}\ge H.
\end{equation}
The above equation is precisely the statement that the level sets of $\rho$ give a sub-inverse-mean-curvature flow. $\qed$

It is relatively easy to see that if a metric is radially monotone everywhere, then it must also be area enlarging everywhere. In particular, the following proposition implies that Kerr-Newman metrics are area enlarging.

\begin{prop}\label{Radially Monotone implies Area Enlarging}
	Let $g$ be an asymptotically flat metric which is everywhere radially monotone. Then $g$ is everywhere area enlarging.
\end{prop}
\textbf{Proof}: Since $g$ is assumed to be globally radially monotone, we have
\begin{equation}
	\frac{\partial(\alpha-2u)}{\partial\rho}\le 0.
\end{equation}
As $g$ is asymptotically flat, we know that
\begin{equation}
	\lim_{\rho\rightarrow\infty}(\alpha-2u)(\rho,z)=0
\end{equation}
for all $z$. Thus, using the fundamental theorem, we may see that
\begin{equation}
	0\le-\int_{\rho_0}^{\infty}\frac{\partial(\alpha-2u)}{\partial\rho}(\rho,z)d\rho=(\alpha-2u)(\rho_0,z).
\end{equation}
This is precisely the coordinate expression of the area enlarging condition. $\qed$
 
We now find several examples of metrics which are area enlarging and strongly asymptotically flat.
\subsection{Axisymmetric Geometrostatic} Here we show that the axisymmetric geometrostatic metrics are area-enlarging and strongly asymptotically flat. Note that the Schwarzschild metric is a member of this family of metrics. Recall that the general form of a geometrostatic metric is
\begin{equation}
	(M,g)=\left(\mathbb{R}^3\backslash\lbrace x_i\rbrace^n_{1},(\chi\psi)^2\delta_{\mathbb{R}^3}\right),
\end{equation}
where for positive numbers $\lbrace a_i\rbrace^{n}_{1}$ and $\lbrace b_i\rbrace^{n}_{1}$ we have
\begin{equation}
	\chi(x)=1+\sum^n_{i=1}\frac{a_i}{\abs{x-x_i}}
\end{equation}
and
\begin{equation}
	\psi(x)=1+\sum^{n}_{i=1}\frac{b_i}{\abs{x-x_i}}.
\end{equation}
If the points $\lbrace x_i\rbrace$ lie on a common line, then the resulting metric will be axisymmetric. The axis of symmetry will be the line on which the $x_i$ lie. After a rotation, we may suppose that the axis of symmetry is the $z$-axis. We may now see that the usual Euclidean cylindrical coordinates are also cylindrical coordinates for $(M,g)$. In particular
\begin{equation}
	g=(\chi\psi)^2(d\rho^2+dz^2+\rho^2d\phi^2).
\end{equation}
A quick calculation shows that the coordinate function $\alpha$ vanishes and
\begin{equation}
	u=-\log(\chi\psi).
\end{equation}
Since both $\chi$ and $\psi$ are strictly larger than one, we see that $u$ is negative. Since $\alpha=0$, it is clear that
\begin{equation}
	\alpha-2u\ge 0.
\end{equation}
This is precisely the coordinate expression of the area-enlarging condition. That $(M,g)$ is also strongly asymptotically flat follows trivially from the fact that $\alpha=0$.

\end{document}